\input amstex
\documentstyle {amsppt}
\NoBlackBoxes

\define \abs#1{\mid #1 \mid}
\define \absq#1{{\mid #1 \mid}^2}
\define \xsing{X_{\text {sing}}}
\define \rn2pi{{ {\sqrt {-1}} \over {2 \pi}}}
\define \ppbar{\partial \overline \partial}
\define \normsq#1{{\mid \mid #1 \mid \mid}^2}
\define \trho{{\tilde \rho}_\alpha}
\define \UP#1{U \times \Bbb P^{#1 - 1}}
\define \FubS{\omega_{\text{Fub-St}}}
\define \OM{\Cal O_{M,m}}
\define \ONf{\Cal O_{N,f(m)}}
\define \fON{f^\prime \Cal O_N}
\define \fONm{(\fON)_m}
\define \fS#1{f^\prime \Cal S_{#1}}
\define \fSm#1{(\fS{#1})_m}
\define \inj#1 #2{0 \rightarrow #1 \rightarrow #2}
\define \tensor #1 #2 #3{{#1} \otimes_{#2} {#3}}
\define \Jx0{\Cal J_{(0,0)}}
\define \OXxpoly{\Cal O_{U,0} [y_0, ... ,y_n]}
\define \OXCx0{\Cal O_{U \times \Bbb C^{n+1},(0,0)}}

\define \frc1N#1{ { {{#1}_1} \over {{#1}_0} }, ... , { {{#1}_n} \over {{#1}_0} } }
\define \Flam{F^{(\lambda)}}

\define \xl0z{(x, \lambda z_0, w)}
\define \Cn1{\Bbb C^{n+1}}
\define \Ctn1{{\tilde {\Bbb C}}^{n+1}}
\define \Pn{\Bbb P^n}
\define \sgst#1{\sigma_{#1}^*}
\define \sinv#1{\sigma_{#1}^{-1}}
\define \zp0{z_0^\prime}

\topmatter

\title
Explicit Construction of Complete K\"ahler Metrics 
of Saper Type by Desingularization
\endtitle

\rightheadtext{Desingularization with One Blow-up \& Saper Metrics}

\author 
Caroline Grant Melles and Pierre Milman
\endauthor

\address 
C. Grant Melles, U. S. Naval Academy, Annapolis, MD 21402-5002
\endaddress

\email cgg\@nadn.navy.mil \endemail

\address 
P. Milman, University of Toronto, Toronto, Ontario  M53 1A1
\endaddress

\email milman\@math.toronto.edu \endemail

\date April 14, 2000  \enddate

\thanks  
The first author was partially supported by a grant from the Naval Academy Research
Council.   The second author was partially supported by 
NSERC grant RGPIN 8949-98 and
the Canada Council Killam Research Fellowship. 
  The authors would also like to
acknowledge support of The Fields Institute for Research in Mathematical Sciences where
a portion of this work was completed. 
\endthanks

\subjclass 32S20, 14E15 \endsubjclass

\abstract  
We construct complete K\"ahler metrics of Saper type  
on the nonsingular set of a subvariety $X$ of a compact K\"ahler manifold 
using 
(a) a method for replacing 
a sequence of blow-ups along smooth centers, used to 
resolve the singularities of $X$, with a single blow-up along a product of coherent 
ideals corresponding to the centers 
and (b) an explicit local 
formula for a Chern form associated to this single blow-up. 
Our metrics have a particularly simple local formula, 
involving essentially a product of distances to the centers of the 
blow-ups used to resolve the singularities of $X$.  
Our proof of (a) uses a generalization of Chow's theorem for 
coherent ideals, proved using the Direct Image Theorem.    
\endabstract
\endtopmatter

\document 

\head Introduction \endhead

Let $X$ be a singular subvariety of a compact K\"ahler manifold $M$.  In 
[GM] we showed how to construct a particular type of complete K\"ahler metric on
the nonsingular set of $X$.  These metrics grow less rapidly than Poincar\'e metrics near the
singular set $\xsing$ of $X$, and are of interest because in certain cases it is known that 
their $L_2$-cohomology equals the intersection cohomology of $X$, while the
$L_2$-cohomology of a Poincar\'e metric is usually not equal to the intersection
cohomology of $X$.  We call our metrics \lq\lq Saper-type" metrics 
after Leslie Saper, who first drew our attention to this subject.  
 Saper proved that on any variety with isolated singularities 
there is a complete K\"ahler metric whose $L_2$-cohomology equals its 
intersection cohomology (see [Sa1], [Sa2]).   Our metrics agree with Saper's in
the case of isolated singularities, but our construction requires no restriction on
the type of singularities.  

The construction of 
Saper-type metrics  in [GM] used  
the geometry of a finite sequence of 
blow-ups along smooth centers 
which resolves the singularities of $X$.      
In this paper we show how to replace a finite sequence of blow-ups along smooth centers 
by a single blow-up along one center (perhaps singular), which we describe in terms 
of its coherent sheaf of ideals $\Cal I$.   Blowing up $M$ along $\Cal I$ desingularizes $X$.  
The support of $\Cal I$ is the singular locus $\xsing$ of $X$.  
We construct $\Cal I$ as a product of coherent ideals $\Cal I_j$  corresponding to the smooth 
centers $C_j$.  Each $\Cal I_j$ is the direct image on $M$ of a product of the 
ideal sheaf of $C_j$ with a sufficiently high power of the exceptional ideal of the 
previous blow-ups.  Our proof that the blow-up of $M$ along $\Cal I$ is equivalent 
to the blow-up along the centers $C_j$ uses a generalization 
of Chow's Theorem for ideals, which we prove using the Direct Image Theorem.  
We then give a simple and explicit construction  of a Chern form 
associated to the blow-up along $\Cal I$, in terms of  local generators of $\Cal I$.  
Finally we use this Chern form to obtain a simpler and more 
explicit form of our Saper-type metrics.  
We also give an example in which we compute $\Cal I$ explicitly in a 
neighborhood of a singular point. 

The Saper-type metric which we obtain can be described in terms of its K\"ahler  
(1,1)-form as 
     $$\omega_S = \omega - \rn2pi \ppbar \log ( \log F)^2,$$
where $\omega$ is the K\"ahler (1,1)-form of a K\"ahler metric on $M$ and $F$ is a
$C^\infty$ function on $M$, vanishing on $\xsing$.  We first construct local $C^\infty$ 
functions $F_\alpha$, on small open sets $U_\alpha$ in $M$, by setting 
     $$F_\alpha = \sum_{j=1}^r \absq {f_j}$$
where $f_1, ... , f_r$ are local holomorphic 
generating functions on $U_\alpha$ for the coherent ideal sheaf $\Cal I$ described above.  
To construct a global metric on $M - \xsing$
(and consequently on $X -\xsing$), we patch with a $C^\infty$ partition of
unity on $M$. 
It is crucial that this patching takes place on $M$, rather than on a
blow-up of $M$, which might add unwanted elements to the $L_2$-cohomology.   

The motivation for our construction is that it may be easier to keep track of a 
single $C^\infty$ function $F$ in a coordinate neighborhood of $p$ in $M$, rather than 
to keep track of the many coordinate neighborhoods associated to successive 
blow-ups which resolve the singularities of $X$.  Similarly, it may be more convenient 
to work with a single ideal sheaf $\Cal I$ on $M$, rather than a sequence of centers and
blow-ups.

\head Table of Contents \endhead

\subhead I.  Outline and Main Results \endsubhead

   \subhead II.  Direct and Inverse Images of Coherent Sheaves of Ideals \endsubhead
   \roster
      \item"{}" Coherent Sheaves
      \item"{}" Direct Images
      \item"{}" Inverse Images
      \item"{}" Composites
      \item"{}" Products of Ideals
   \endroster

   \subhead III.  Blowing up a Complex Manifold along a Coherent Sheaf of Ideals \endsubhead
   \roster
      \item"{}" Local Description of Blow-ups
      \item"{}" Global Description of Blow-ups
     \item"{}"  Ideals, Divisors, Line Bundles, and Sections
     \item"{}"  Exceptional Divisors of Blow-ups
     \item"{}"  Exceptional Line Bundles of Blow-ups
     \item"{}"  Universal Property of Blow-ups
     \item"{}"  Blow-up of a Product of Ideals
     \item"{}"  Direct Images under Blow-up Maps  
 \endroster

   \subhead IV.  Chow's Theorem for Ideals \endsubhead

   \subhead V.  Chow's Theorem Applied to Blow-ups  \endsubhead 
    
  \subhead VI.  Replacing a Sequence of Blow-ups by a Single Blow-up \endsubhead

   \subhead VII.  Chern Forms and Metrics for Exceptional Line Bundles \endsubhead
 \roster
   \item"{}" Chern Forms on Line Bundles  
   \item"{}"  Local Chern Forms for Blow-ups
   \item"{}" Global Chern Forms for Blow-ups
\endroster

   \subhead VIII.  Construction of Saper-Type Metrics \endsubhead
\roster
   \item"{}"  Local Construction of Metrics
   \item"{}"  Global Construction of Metrics
\endroster

\subhead IX.  Example \endsubhead

\head I.  Outline and Main Results \endhead

In sections II and III we give some background and basic results about coherent sheaves 
of ideals and blow-ups.  
We begin by describing the direct and inverse images of  sheaves, and 
in particular, direct and inverse images of coherent sheaves of ideals. Then 
we describe the blow-up $\pi : \tilde M \rightarrow M$ of a complex manifold $M$ along a 
coherent sheaf of ideals $\Cal I$.  The analytic subset $C = V(\Cal I)$ of $M$ 
determined by $\Cal I$ is 
called the center of the blow-up.   If $C$ is smooth and of codimension 
at least 2, then $\tilde M$ is smooth.  
The blow-up map $\pi$ is proper and 
is a  biholomorphism except along its exceptional divisor $E = \pi^{-1} (C)$. 
Even though the direct image of an ideal sheaf  
may not be an ideal sheaf in general, 
the direct image of an ideal sheaf  under a blow-up map {\bf is} an 
ideal sheaf.     

Section IV is devoted to a proof of Chow's Theorem for Ideals using 
the Direct Image Theorem, which states that the direct image of a coherent sheaf under 
a proper map is coherent.
Section V 
contains some corollaries for blow-up maps which are useful in constructing 
single-step blow-ups from a sequence of blow-up maps.  
 
\proclaim{Chow's Theorem for Ideals}  Let $U$ be an open neighborhood of $0$ 
in $\Bbb C^m$ and let $X$ be an analytic subset of $U \times \Bbb P^n$.  Let 
$\Cal J$ be a coherent sheaf of ideals on $X$.  Then $\Cal J$ is {\bf relatively 
algebraic} in the following sense:  $\Cal J$ is generated (after shrinking $U$ 
if necessary) by a finite number of homogeneous polynomials in homogeneous 
$\Bbb P^n$-coordinates, with analytic coefficients in $U$-coordinates.
\endproclaim

Chow's Theorem for Ideals helps to describe the relatively algebraic
structure of blow-ups.   
The most useful corollary for the purposes of this paper is the following, which shows 
that, even though the inverse image of the direct image of an ideal sheaf may not 
be the original ideal sheaf in general, on a blow-up of a compact complex manifold 
we can ensure that the two are equal by first multiplying by a high enough 
power of the ideal sheaf $\Cal I_E$ of the exceptional divisor.  
This corollary is similar to results of Hironaka and Rossi in [HR] but 
our proof uses simpler and more explicit methods and is more constructive 
in nature.   We go on to apply this corollary repeatedly to get 
an explicit description of a coherent sheaf for single-step blow-ups, 
as a product of coherent sheaves corresponding to a sequence of 
blow-ups along smooth centers.

\proclaim{Corollary}   Let $\pi : \tilde M \rightarrow M$ be the blow-up of a 
compact complex manifold $M$ along a coherent sheaf of ideals  $\Cal J_1$ 
and let $E$ be the exceptional divisor of $\pi$.  Let $\Cal J_2$ be a coherent sheaf 
of ideals on $\tilde M$.  Then there exists an integer $d_0$ such that 
     $$\pi^{-1} \pi_* ( \Cal J_2 \Cal I_E^d ) = \Cal J_2 \Cal I_E^d$$ 
for all $d \geq d_0$.  
\endproclaim

For the purposes of this paper and to apply Hironaka's theorem on 
embedded resolution of singularities, we need only consider blow-ups 
of smooth spaces.  
If $\tilde M$ is smooth, 
the blow-up of $\tilde M$ along $\Cal J_2 $ is isomorphic to 
the blow-up of $\tilde M$ along $\Cal J_2 \Cal I_E^d$.  Furthermore, 
the blow-up of $\tilde M$ along 
$\Cal J_2 \Cal I_E^d $ is isomorphic to the blow-up of the base space $ M$ 
along $\Cal J_1 \pi_*(\Cal J_2 \Cal I_E^d) $.  Thus we can replace the 
pair of blow-ups, first along $\Cal J_1$ and then along $\Cal J_2$,  
by a single blow-up along $ \Cal J_1 \pi_* (\Cal J_2 \Cal I_E^d)$.  
Repeating this procedure for a finite sequence of 
smooth centers enables us to construct a coherent sheaf of ideals $\Cal I$
on $M$ such blowing up $M$ along $\Cal I$ is equivalent to 
blowing up successively along smooth centers.  Section VI contains a more detailed version of the proof of the following proposition.

\proclaim{Proposition (Single-Step Blow-ups)}  Let $M$ be a compact complex manifold 
and let 
     $$M_m \overset {\pi_m} \to \rightarrow M_{m-1} 
     \rightarrow ... \rightarrow
     M_2 \overset {\pi_2} \to \rightarrow M_1 
     \overset {\pi_1} \to \rightarrow M_0 = M$$ 
be a finite sequence of blow-ups along smooth centers $C_j \subset M_{j-1}$ 
of codimension at least 2.  
Then there is a coherent sheaf of ideals $\Cal I$ on $M$ such that the 
blow-up of $M$ along $\Cal I$ is isomorphic to the blow-up   
of $M$ along the sequence of smooth centers $C_j$.  Furthermore, 
we may construct $\Cal I$ to be of the form  
     $$\Cal I = \Cal I_1 \Cal I_2 ... \Cal I_m,$$
where each $\Cal I_j$ is a coherent sheaf of ideals on $M$ and   
\roster 
   \item"{i.}"  $\Cal I_j$ is the direct image on $M$ of the 
ideal sheaf of $C_j$ multiplied by a high enough power of 
the ideal sheaf of the exceptional divisor of the first $j-1$ blow-ups, 
\item"{ii.}"  
the inverse image of $\Cal I_j$ on $M_{j-1}$ is 
the ideal sheaf of $C_j$ multiplied by the same power of the 
exceptional ideal sheaf as in (i), and 
\item"{iii.}"  the blow-up of $M_{j-1}$ along the inverse 
image of $\Cal I_j$ is isomorphic to the blow-up of $M_{j-1}$ 
along $C_j$. 
\endroster
\endproclaim

This result is related to Theorem II.7.17 of  [Ha1], but our 
proposition is much 
more explicit and constructive in nature.  

 We are particularly interested in the case of a sequence of blow-ups 
along smooth centers 
which resolves the singularities of a singular subvariety $X$ of $M$.   
In this case, the proposition gives us a coherent ideal sheaf $\Cal I$ 
on $M$, supported on the singular locus of $X$, such that blowing up 
along $\Cal I$ desingularizes $X$, and also gives a factorization of $\Cal I$ 
in terms of ideals corresponding to the original sequence of blow-ups.  
This factorization of $\Cal I$ is essentially unique for curves.

In section VII we give a simple and explicit construction of a Chern form 
associated to a blow-up.  
Suppose that  $\pi : \tilde M \rightarrow M$ is  
the blow-up of a complex manifold $M$ along a coherent sheaf of ideals  
$\Cal I$ such that $\tilde M$ is smooth.  
Let $E$ be the exceptional divisor and $L_E$ the line 
bundle on $\tilde M$ associated to $E$.  Let $f_1, ..., f_r$ be local 
holomorphic generating functions for $\Cal I$ on a small open set $U \subset M$. 
We construct a Chern form for $L_E$ on $\tilde U = \pi^{-1} (U)$ 
by pulling back the negative of a Fubini-Study form on projective space.  
This Chern form is strictly negative on the fibres of the map 
$E \rightarrow C = V(\Cal I)$, and is given on $\tilde U - \tilde U \cap E$ 
by 
     $$c_1 (L_E) = \pi^* (- \rn2pi  \ppbar \log \sum_{j=1}^r \absq {f_j(z)}).$$
If $M$ is compact, we may patch together local Chern forms using a $C^\infty$ 
partition of unity on $M$, in such a way that the negativity on fibres is preserved.  

Now consider in more detail 
a singular subvariety $X$ of a compact K\"ahler manifold $M$.  
Hironaka's famous theorem  tells us that the singularities of $X$ may 
be resolved by a finite sequence of blow-ups of $M$ along smooth centers, 
such that the total exceptional divisor of the composite of all the blow-ups 
is a normal crossings divisor $D$ in $\tilde M$ which has normal crossings with 
the desingularization $\tilde X $ in $\tilde M$ and such that 
$\tilde M - D \cong M - \xsing$ and $\tilde X - 
\tilde X \cap D \cong X - \xsing$ (cf. [BM]).    
By the Single-Step Blow-up Proposition, we may resolve the singularities 
of $X$ by blowing up $M$ along a single coherent sheaf of ideals $\Cal I$ on 
$M$, whose blow-up is isomorphic to the blow-up obtained using the sequence 
of smooth centers.  The inverse image ideal sheaf of $\Cal I$ in the blow-up $\tilde M$ 
determines the normal crossings divisor $D$ and the support of $\Cal I$ in 
$M$ is $\xsing$.  
We construct a Chern form for the blow-up along $\Cal I$, using local holomorphic 
generating functions of $\Cal I$ as above and patching with a $C^\infty$ 
partition of unity on $M$.  This Chern form is negative definite on the 
fibres of the map from $D$ to $\xsing$.  Subtracting this Chern form from 
the K\"ahler (1,1)-form 
of a K\"ahler metric on $M$ gives  
the (1,1)-form of a K\"ahler metric on $\tilde M$, our \lq\lq desingularizing  
metric."  The completion of $X - \xsing$ with respect to this metric is 
nonsingular.  Similarly, we  use the local holomorphic generators of $\Cal I$ 
to construct our complete K\"ahler Saper metric on $\tilde M - D \cong 
M - \xsing$.  Both metrics are described in more 
detail below and in section VIII.

\proclaim{Theorem}  
Let $X$ be a singular subvariety of a compact K\"ahler manifold $M$ and let $\omega$
be the K\"ahler (1,1)-form of a K\"ahler metric on $M$.  Then there exists a $C^\infty$ function
$F$ on $M$,
vanishing on $\xsing$, such that for $l$ a large enough positive integer, 
 \roster
   \item"{i.}" the (1,1)-form 
     $$\tilde \omega = l \omega + \rn2pi \ppbar \log F$$ 
is the K\"ahler form of a desingularizing K\"ahler metric for $X$, i.e.
the completion of $X - \xsing$ with respect to $\tilde \omega$ is a desingularization of
$X$ and 
   \item"{ii.}" the (1,1)-form 
     $$\omega_S = l \omega - \rn2pi \ppbar \log ( \log F)^2$$ 
is the K\"ahler form of a complete K\"ahler modified Saper metric (in the
terminology of [GM]) on $M - \xsing$ and hence on $X - \xsing$.
\endroster
Furthermore, the function $F$ may be constructed 
to be of the form 
     $$F = \prod_\alpha F_\alpha^{\rho_\alpha},$$ 
where $\{ \rho_\alpha \}$ is a $C^\infty$ partition of 
unity subordinate to an open cover $\{U_\alpha \}$ of $M$, 
$F_\alpha$ is a function on $U_\alpha$ of the form 
     $$F_\alpha = \sum_{j=1}^r \absq {f_j},$$ 
and $f_1, ... , f_r$ are holomorphic functions on $U_\alpha$, 
vanishing exactly on $\xsing \cap U_\alpha$.  
More specifically, $f_1, ... , f_r$ are local holomorphic 
generators of a coherent ideal sheaf $\Cal I$ on $M$ 
such that blowing up $M$ along 
$\Cal I$ desingularizes $X$, 
$\Cal I$ is supported on $\xsing$, 
and the exceptional 
divisor of the blow-up along $\Cal I$ has normal crossings with itself and 
with the desingularization of $X$.     
\endproclaim

The coherent ideal sheaf $\Cal I$ is constructed as a product 
$\Cal I_1 \Cal I_2 ... \Cal I_m$ of coherent ideal sheaves 
corresponding to a sequence of blow-ups along smooth centers $C_j$ 
which resolves the singularities of $X$.  This factorization 
of $\Cal I$  
gives  a corresponding factorization of $F_\alpha$, as essentially a product of 
distances to the centers,    
     $$ F_\alpha = \prod_{j=1}^m \sum_{i=1}^{r_j} \absq {v_{j i}}$$
where, for each $j$, the functions $\{v_{j i}\}$ are 
local holomorphic functions on $U_\alpha$ whose
pullbacks to the preimage of $U_\alpha$  
under the first $j-1$ blow-ups 
generate an ideal sheaf with the same blow-up as
$C_j$.  

The idea behind the  metric constructions in this paper is to first find simple and explicit
formulas locally on $M$, then patch by $C^\infty$ partitions of unity on $M$.  We wish
to avoid formulas which are local only on blow-ups  of $M$ and we also wish to
avoid introducing $C^\infty$ partition-of-unity functions on the blow-ups.  
                                                  
     We conclude, in section IX, by constructing $\Cal I$ for the cuspidal cubic $y^2 - x^3$.  
The method used generalizes to the case of 
any singular locally toric complex analytic variety. 
The details 
will be given elsewhere.  

\head II.  Direct and Inverse Images of Coherent Sheaves of Ideals  \endhead

\head  Coherent Sheaves  \endhead

We first review the important concept of coherence (see e.g. [GrR1], [GuR]). 

Let $M$ be a complex space and let $\Cal S$ be an analytic sheaf on $M$, i.e. a sheaf of
$\Cal O_M$-modules. For example, consider an ideal sheaf of ${\Cal O}_M$ or the sheaf of
holomorphic
sections of a holomorphic vector bundle on $M$. 

\definition{Definition}  
The sheaf $\Cal S$ is {\bf of finite type at $x \in M$} if there exists an open set $U$ of
$x$ such that the restriction $\Cal S_U$ of $\Cal S$ to $U$ is generated by a finite
number of sections of $\Cal S$ over $U$.  This means that there exist sections $s_1, ... ,
s_r$ of $\Cal S$ over $U$ such that for each point $y \in U$ and for each germ $g_y \in
\Cal S_y$, there exist $a_{1y}, ... ,a_{ry} \in \Cal O_{M,y}$ such that
     $$g_y = \sum_{i=1}^r a_{iy} s_{iy}.$$
The sheaf $\Cal S$ is {\bf of finite type on $M$} if $\Cal S$ is of finite type at $x$ for all $x \in
M$. 
\enddefinition

\remark{Remark}  
Note that if $s$ and $t$ are sections of $\Cal S$ on a neighborhood of a point $y$ 
such that $s_y = t_y$ (i.e. they have the same germs at $y$), 
then $s= t$ in an open neighborhood of $y$, by fundamental properties of sheaves.  In
particular, in the definition above, if $g_y, a_{1y}, ... , a_{ry}$ are the germs of $g, a_1,
... ,a_r$ at $y$ then there exists a neighborhood $V \subset U$ of $y$ such that
     $$g = \sum_{i=1}^r a_i s_i$$ on $V$.
\endremark

Each finite collection $s = (s_1, ... , s_r)$ of sections of $\Cal S$ over $U$ determines a
sheaf homomorphism
     $$\psi_s : \Cal O^r_U \rightarrow \Cal S_U$$
given by
     $$(f_1, ... , f_r) \mapsto \sum_{i=1}^r f_i s_i.$$

\definition{Definition} The sheaf 
$\Cal S$ is {\bf of relation finite type at $x \in M$} if ker $\psi_s$ 
is of finite type at
$x$ for all finite collections $s$ of sections of $\Cal S$ over an open neighborhood $U$
of $x$.  $\Cal S$ is of {\bf relation finite type on $M$} if $\Cal S$ is of relation finite
type at $x$ for all $x \in M$. 
\enddefinition

\definition{Definition} The sheaf 
$\Cal S$ is {\bf coherent on $M$} if
 \roster
   \item  $\Cal S$ is of finite type on $M$, and
   \item  $\Cal S$ is of relation finite type on $M$.
\endroster
\enddefinition

Since coherent sheaves are always finite type, by definition, it follows that 
if $\Cal S$ is a coherent sheaf on
a complex space $X$ and $s_1, ... ,s_r$ are sections of $\Cal S$ on a neighborhood $U$
of a point $x$ such that the germs $s_{1x}, ... , s_{rx}$ generate $\Cal S_x$, then there
exists a neighborhood $V \subset U$ of $x$ such that $s_1, ... ,s_r$ generate $\Cal
S_V$.

We refer the reader to [F], [GrR1], [GrR2], [GuR], and [W]  
for background on the following and other fundamental
properties of coherent sheaves: 
\roster 
   \item"{i.}"  The sheaf $\Cal O_M$ is coherent.  
   \item"{ii.}"  A subsheaf of a coherent sheaf is coherent 
if and only if it is of finite type.  In
particular, an ideal sheaf of $\Cal O_M$ is coherent if and only if 
it is of finite type.
     \item"{iii.}" A coherent ideal sheaf $\Cal I$ on a complex space determines a 
closed complex analytic subspace $V(\Cal I)$, and the ideal sheaf $\Cal I_Y$ of a 
closed complex analytic subspace $Y$ of a complex space is coherent. 
\endroster

\proclaim{Lemma II.1}  If $\Cal I_1$ and $\Cal I_2$ are coherent sheaves of 
ideals on a complex space $M$, then the product ideal sheaf $\Cal I_1 \Cal I_2$ is 
also coherent.
\endproclaim

\demo{Proof}  Since both $\Cal I_1$ and $\Cal I_2$ are of finite type, their 
product is of finite type and is thus coherent. \qed
\enddemo

We define direct images and inverse images of
coherent sheaves of ideals, and give some conditions under which these sheaves are themselves
coherent ideal sheaves (in general they may be only sheaves of modules).  We 
show that direct and inverse images of composite maps are composites of the direct  
and inverse image maps (functoriality).  We also show that the inverse image of a product 
of ideals is the product of the inverse image ideals.  
Direct and inverse images of ideal sheaves under blow-up maps are discussed 
in Lemmas III.9 and V.8.

\head Direct Images \endhead

\definition{Direct Image}  
Let $f: M \rightarrow N$ be a holomorphic map of complex spaces and let $\Cal S$
be a sheaf on $M$.  The direct image sheaf $f_* \Cal S$ on $N$ is the sheaf associated
with the presheaf 
given by $f_* \Cal S (U) = \Cal S ( f^{-1} (U))$, for $U$ any open set in $N$. 
\enddefinition

If $\Cal S$ is coherent, 
the direct image $f_* \Cal S$ is not necessarily coherent.  However $f_* \Cal S$ is
coherent if $f$ is proper, by the Direct Image Theorem.    
We recall the Direct Image Theorem in our context (see e.g. [GrR1], pp 207, 227,  and 36). 

\proclaim{Direct Image Theorem}  
Let $f : M \rightarrow N$ be a holomorphic map of complex spaces and let $\Cal S$
be a coherent sheaf on $M$.  If $f$ is proper then $f_* \Cal S$ is coherent. 
\endproclaim

In particular, 
if $f$ is a blow-up map (see section III), then $f$ is proper and $f_* \Cal S$ is 
coherent if $\Cal S$ is.

If $\Cal J$ is a sheaf of ideals on $M$, then $f_* \Cal J$ is a sheaf of $\Cal
O_N$-modules but not, in general, an ideal sheaf on $N$.  We will show  
(Lemma III.9)  that if $f$ is a blow-up map then $f_* \Cal J$ is an ideal sheaf.  

\head Inverse Images \endhead
     
Once again, let $f: M \rightarrow N$ be a holomorphic map of complex spaces.  Let
$\Cal S$ be a sheaf of $\Cal O_N$-modules. 

\definition{Topological Inverse Image} 
We define the topological inverse image $f^\prime \Cal S$ to be the fibre product $\Cal S
\times_N M$, i.e. the stalk of $f^\prime \Cal S$ over a point $m \in M$ is the stalk of
$\Cal S$ over $f(m) \in N$:
     $$(f^\prime \Cal S)_m = \Cal S_{f(m)}.$$
\enddefinition

Note that $f^\prime \Cal S$ is a sheaf of $f^\prime \Cal O_N$-modules.  If $\Cal S$ is
coherent then so is $f^\prime \Cal S$. 

\definition{Pullback Sheaf}  
We define the pullback sheaf as 
     $$f^* \Cal S 
     = f^\prime \Cal S \otimes_{f^\prime \Cal O_N} \Cal O_M.$$
\enddefinition

Note that $f^* \Cal S$ is a sheaf of $\Cal O_M$-modules and once again, if $\Cal S$ is
coherent then so is $f^* \Cal S$.  Also
     $$f^* \Cal O_N = \tensor {\fON} {\fON} {\Cal O_M} = \Cal O_M.$$

If $\Cal I$ is an ideal sheaf on $N$, we have an exact sequence
     $$ 0 \rightarrow \Cal I \rightarrow \Cal O_N.$$
Since tensoring is not in general left exact, the induced map 
     $$f^* \Cal I \rightarrow f^* \Cal O_N =\Cal O_M$$ 
is not necessarily injective, so $f^* \Cal I$ is not necessarily an {\it ideal sheaf} on $M$. 
The {\it image} of $f^* \Cal I$ in $\Cal O_M$ is an ideal sheaf, which we call the 
inverse image ideal sheaf and will describe in more detail later in this section.  

\subhead Flat Maps \endsubhead   A holomorphic map $f: M \rightarrow N$ of complex
spaces is  flat if

       $$\OM \text{ is $\ONf$-flat}$$ 
for all $m \in M$.  Equivalently, $f$ is flat if 
for every exact sequence
     $$\inj {S_1} {S_2}$$
 of $\ONf$-modules, the
induced sequence 
     $$\inj {\tensor {S_1} {\ONf} {\OM} } {\tensor {S_2} {\ONf} {\OM}}$$ 
is also exact. 

     There are many references on flat maps, e.g. ([F], p. 147 and p. 155).

\example{Examples}  
If $X$ and $Y$ are complex spaces, the canonical projection $X \times Y \rightarrow
Y$ is flat.  Every locally trivial holomorphic map is flat.  
In particular, if $f :  L \rightarrow X$ is a line bundle over a complex
space $X$ (or more generally, a vector bundle), then $f$ is flat. 
\endexample

\proclaim{Lemma II.2}  
If $f: M \rightarrow N$ is 
a flat holomorphic map of complex spaces 
and $\inj {\Cal S_1} {\Cal S_2}$ is an exact sequence of
sheaves of $\Cal O_N$-modules, then $\inj {f^* \Cal S_1} {f^* \Cal S_2}$ is an exact
sequence of sheaves of $\Cal O_M$-modules. 
\endproclaim

\demo{Proof}  
Suppose that 
     $$\inj {\Cal S_1} {\Cal S_2}$$
is an exact sequence of sheaves of $\Cal O_N$-modules, i.e. 
     $$\inj {\Cal S_{1,n}} {\Cal S_{2,n}}$$
is an exact sequence of $\Cal O_{N,n}$-modules for each $n \in N$.  
Then in particular, 
     $$\inj {\Cal S_{1,f(m)}} {\Cal S_{2,f(m)}}$$
is an exact sequence of $\ONf$-modules for all $m \in M$. 
If $f : M \rightarrow N$ is flat, then
     $$\inj {\tensor {S_{1,f(m)}} {\ONf} {\OM} } {\tensor {S_{2,f(m)}} {\ONf} {\OM}}$$ 
is exact for all $m \in M$, i.e. 
     $$ \inj {\tensor {\fSm 1} {\fONm} {\OM}}
     {\tensor {\fSm 2} {\fONm} {\OM}}$$
is exact for all $m \in M$.  These tensor products can be rewritten as
     $$ \inj {(\tensor {\fS 1} {\fON} {\Cal O_M})_m}
     ({\tensor {\fS 2} {\fON} {\Cal O_M})_m},$$
showing that 
     $$\inj {\tensor {\fS 1} {\fON} {\Cal O_M}}
     {\tensor {\fS 2} { \fON} {\Cal O_M}}$$
is exact.  By the definition of $f^*$, this means that 
     $$\inj {f^* \Cal S_1} {f^* \Cal S_2}$$
is exact.
\qed  
\enddemo

\proclaim{Lemma II.3}  
If $\Cal L$ is the sheaf of holomorphic sections of a line bundle (or more generally of a
vector bundle) over a complex space $M$, and
     $$\inj {\Cal S_1} {\Cal S_2}$$
is an exact sequence of sheaves of $\Cal O_M$-modules, then
     $$\inj {\Cal S_1 \otimes \Cal L} {\Cal S_2 \otimes \Cal L}$$ 
is also exact.
\endproclaim

\demo{Proof}  
A finitely generated module over a local noetherian ring is flat if and only if it is free
([Ma],  Proposition 3.G,  p. 21).  Therefore $\otimes_{\OM} \Cal L_m$ preserves
exact sequences. 
\qed
\enddemo

\definition{Inverse Image Ideal}  
Let $f: M \rightarrow N$ be a holomorphic map of complex spaces.  
If $\Cal I$ is an ideal sheaf on $N$, the image of $f^* \Cal I$ in $\Cal O_M$ is 
an ideal sheaf which we define to be the inverse image ideal sheaf $f^{-1} \Cal I$.    
\enddefinition

The ideal sheaf $f^{-1} \Cal I$ is sometimes written $f^* \Cal I \cdot \Cal O_M$ 
or $f^{-1} \Cal I \cdot \Cal O_M$.  If
$\Cal I$ is coherent, then $f^{-1} \Cal I$ is also coherent. 

If $\Cal I$ is a coherent ideal, the subscheme of $M$ determined by $f^{-1} \Cal I$ is the
inverse image scheme of the subscheme of $N$ determined by $\Cal I$, i.e. 
     $$V(f^{-1} \Cal I) = f^{-1} (V(\Cal I)).$$

\proclaim{Lemma II.4}  
If $f: M \rightarrow N$ is a flat holomorphic map of complex spaces 
and $\Cal I$ is an ideal sheaf on $N$, 
then $f^{-1} \Cal I \cong f^* \Cal I$.
\endproclaim

\demo{Proof}  
By Lemma II.2 above, if $f$ is flat, then the map $f^* \Cal I \rightarrow f^* \Cal O_N =
\Cal O_M$ is injective. 
\qed 
\enddemo

\proclaim{Corollary II.5}  
If $f: L \rightarrow X$ is a line bundle (or more generally a vector bundle) 
and $\Cal I$ is an ideal sheaf on $X$, then $f^{-1}
\Cal I = f^* \Cal I$. 
\endproclaim

\demo{Proof}  
As noted in the discussion of flat maps above, the projection of a line bundle (or vector
bundle) onto its base space is a flat map. 
\qed  
\enddemo

\head Composites \endhead

Next we describe the behavior of direct and inverse images under composites.  The 
proofs are straightforward, using the definitions above.  

\proclaim{Lemma II.6 (The Composite of Direct Images is the Direct Image 
of the Composite)} Let $M_1 \overset f \to \rightarrow M_2 
\overset g \to \rightarrow M_3$ be holomorphic maps of complex spaces and 
let $\Cal S$ be a sheaf on $M_1$.  Then 
     $$g_* (f_* \Cal S) \cong (g \circ f)_* \Cal S.$$
\endproclaim 

\demo{Proof} Let $U$ be an open set in $M_3$.  Then 
     $$\align   g_* (f_* \Cal S) (U) & = 
       (f_* \Cal S) (g^{-1} (U)) \\
       & = \Cal S (f^{-1} g^{-1} (U)) \\
       & = \Cal S ((g \circ f)^{-1} (U)) \\
       & = (g \circ f)_* (U). \qed \endalign$$
\enddemo

\proclaim{Lemma II.7 (The Composite of Topological Inverse Images is 
the Topological Inverse Image of the Composite)} 
Let $M_1 \overset f \to \rightarrow M_2 
\overset g \to\rightarrow M_3$ be holomorphic maps of complex spaces and 
let $\Cal S$ be a sheaf on $M_3$.  Then 
     $$f^{\prime} (g^{\prime} \Cal S) \cong (g \circ f)^{\prime} \Cal S.$$
\endproclaim

\demo{Proof} We will prove the statement on stalks.  Let $m$ be a point in 
$M_1$.   Then
     $$\align (f^\prime (g^\prime \Cal S))_m & = (g^\prime \Cal S)_{f(m)} \\
       & = {\Cal S}_{g \circ f (m) } \\
       & = ( (g \circ f)^\prime \Cal S)_m. \qed \endalign $$
\enddemo

\proclaim{Lemma II.8 (The Composite of Pullbacks is the Pullback of the 
Composite)}  Let $M_1 \overset f \to \rightarrow M_2 
\overset g \to \rightarrow M_3$ be holomorphic maps of complex spaces and 
let $\Cal S$ be a sheaf on $M_3$.  Then 
     $$f^{*} (g^* \Cal S) \cong (g \circ f)^{*} \Cal S.$$
\endproclaim 

\demo{Proof} For convenience, let $\Cal O_i$ represent $\Cal O_{M_i}$ for $i = 1,2,3$.  
Recall that 
       $$g^* \Cal S = g^\prime \Cal S \otimes_{g^\prime \Cal O_3} \Cal O_2.$$
Similarly  
     $$\align f^* (g^* \Cal S) 
       & = f^\prime (g^* \Cal S) \otimes_{f^\prime \Cal O_2} \Cal O_1 \\
       & = f^\prime (g^\prime \Cal S \otimes_{g^\prime \Cal O_3} \Cal O_2) 
         \otimes_{f^\prime \Cal O_2} \Cal O_1. \endalign$$
Looking at stalks over $m \in M_1$ we have 
     $$\align (f^*(g^* \Cal S))_m 
     & =  f^\prime (g^\prime \Cal S \otimes_{g^\prime \Cal O_3} \Cal O_2)_m 
         \otimes_{(f^\prime \Cal O_2)_m} \Cal O_{1,m}\\
       & = (g^\prime \Cal S \otimes_{g^\prime \Cal O_3} \Cal O_2)_{f(m)} 
       \otimes_{\Cal O_{2,f(m)}} \Cal O_{1,m} \\
       & = (g^\prime \Cal S)_{f(m)} \otimes_{(g^\prime \Cal O_3)_{f(m)}} 
       \Cal O_{2, f(m)}  \otimes_{\Cal O_{2,f(m)}} \Cal O_{1,m} \\
       & = \Cal S_{g(f(m))} \otimes_{\Cal O_{3,g(f(m))} } 
       \Cal O_{2, f(m)}  \otimes_{\Cal O_{2,f(m)}} \Cal O_{1,m} \\
       & = \Cal S_{g(f(m))} \otimes_{\Cal O_{3,g(f(m))} } \Cal O_{1,m} \\
       & = ((g \circ f)^\prime \Cal S)_m \otimes_{((g \circ f)^\prime \Cal O_3)_m} 
       \Cal O_{1,m} \\
     & = ((g \circ f)^* \Cal S)_m. \qed \endalign$$
\enddemo

The following lemma is more naturally understood in terms of subschemes 
determined by coherent sheaves of ideals.  Its interpretation in term of 
subschemes is that the inverse image subscheme under a composite map is 
the composite of the inverse images.  Briefly, $f^{-1}$ is functorial  
on ideals and their corresponding subschemes.  

\proclaim{Lemma II.9 (The Composite of Inverse Images is the Inverse Image of the 
Composite)}  Let $M_1 \overset f \to \rightarrow M_2 
\overset g \to \rightarrow M_3$ be holomorphic maps of complex spaces and 
let $\Cal I$ be a sheaf of ideals on $M_3$.  Then 
     $$f^{-1} (g^{-1} \Cal I) \cong (g \circ f)^{-1} \Cal I.$$
\endproclaim 
       
\demo{Proof} As in the previous proof, let $\Cal O_i = \Cal O_{M_i}$.  
Recall that $g^{-1} \Cal I$ is defined to be the image of $g^* \Cal I$ in 
$\Cal O_2$, so there is a surjective map 
     $$g^* \Cal I \mapsto g^{-1} \Cal I.$$
The map of topological inverse images 
     $$f^\prime g^* \Cal I \mapsto f^\prime g^{-1} \Cal I$$ 
is also surjective.  

     Tensoring over $f^\prime \Cal O_2$ by $\Cal O_1$ we obtain the map 
     $$f^* g^* \Cal I \mapsto f^* g^{-1} \Cal I,$$
which is surjective since tensoring is right exact.  

     Finally we note that 
     $$\alignat2  f^{-1} g^{-1} \Cal I 
       &= \text{image of $f^* g^{-1} \Cal I$ in $\Cal O_1$} 
       &&\qquad \text{by definition} \\
       &= \text{image of $f^* g^* \Cal I$ in $\Cal O_1$} 
       &&\qquad \text{by surjectivity} \\
       &= \text{image of $(g \circ f)^* \Cal I$ in $\Cal O_1$}
       &&\qquad \text{by Lemma II.8}\\
       &= (g \circ f)^{-1} \Cal I
       &&\qquad \text{by definition}  \qed \endalignat$$
       
\enddemo

\head Products of Ideals \endhead

The following lemma is also more naturally understood in terms of subschemes 
determined by coherent sheaves of ideals. 
The subscheme of $M$ determined by $(f^{-1} \Cal I_1) (f^{-1} \Cal I_2)$ 
is the union of the subschemes determined by $f^{-1} \Cal I_1$ and $f^{-1} \Cal  I_2$, 
which are the inverse images of the subschemes determined by $\Cal I_1$ and  $\Cal I_2$.   
The subscheme of $M$ determined by $f^{-1} (\Cal I_1 \Cal I_2)$ is the inverse image 
of the union of the subschemes determined by $\Cal I_1$ and $\Cal I_2$, which 
is the same as the union of the inverse images.  
     
\proclaim{Lemma II.10 (The Inverse Image Ideal of a Product of Ideal Sheaves is the 
Product of the Inverse Image Ideal Sheaves)} Let $f: M  \rightarrow 
 N$ be a holomorphic 
map of complex spaces and let $\Cal I_1$ and $\Cal I_2$ be sheaves of ideals on $N$.  
Then
     $$(f^{-1} \Cal I_1) (f^{-1} \Cal I_2) \cong f^{-1} (\Cal I_1 \Cal I_2).$$
\endproclaim

\demo{Proof}  Note that both $f^{-1} (\Cal I_1 \Cal I_2)$ and $(f^{-1} \Cal I_1) 
(f^{-1} \Cal I_2)$ are generated as ideals in $\Cal O_M$ by products of 
the form $f^* w_1 f^* w_2$ where $w_1$ is a germ of $\Cal I_1$ and $w_2$ a germ 
of $\Cal I_2$. 
\qed
\enddemo

The direct image of a product of ideal sheaves is not necessarily equal to the product 
of the direct images, but we will show later (Lemma V.8) that the two are equal if the map 
is a blow-up of a smooth center and the ideal sheaves are first multiplied by a high 
enough power of the ideal sheaf of the exceptional divisor.

\head III.  Blowing up a Complex Manifold along a Coherent Sheaf of Ideals \endhead

Let $M$ be a complex manifold and let $\Cal I$ be a coherent sheaf of ideals 
on $M$.  Here and throughout the paper we will always assume that $\Cal I$ 
is not the zero sheaf.  Since $\Cal I$ is coherent, 
for each point $p \in M$ we may
choose a coordinate neighborhood $U$, centered at $p$, such that $\Cal I(U)$ is
generated by a finite number of global sections over $U$.  
We first define the blow-up of $M$ along $\Cal I$ locally over such an open set
$U$, using a collection of generators of $\Cal I(U)$.  We then show that the result is
independent of the collection of generators chosen, so that the blow-up may be defined
globally over $M$.  

Blow-ups may also be defined for singular complex spaces but we do not need such
generality here. 

\head  Local Description of Blow-ups  \endhead

Let $M$ be a complex manifold and $\Cal I$ a coherent sheaf of ideals on $M$ as above. 
Let $U$ be a small enough coordinate neighborhood 
in $M$ that $I = \Cal I (U)$ is generated by a finite collection of 
global sections $f_1, ... , f_r$ on $U$.  
Set
     $$V(I) = \{ z \in U  : h(z) =0 \quad \text{for all} \quad  h \in I \}.$$
We define a map
      $$\phi_f : U - V(I) \rightarrow \Bbb P^{r-1}$$ 
by setting $\phi_f (z) = [f_1 (z): ... : f_r (z)]$.  Let $\Gamma (\phi_f)$ be the graph of
$\phi_f$ in $U \times \Bbb P^{r-1}$, i.e. 
     $$\align \Gamma (\phi_f) 
     &= \{ (z, [\xi]) : z \in U - V(I) 
     \quad \text{and} \quad [\xi] = [f_1 (z): ... :f_r(z)]\}  \\
     &= \{ (z, [\xi]): z \in U - V(I) 
     \quad \text{and} \quad 
     f_i(z) \xi_j = f_j (z) \xi_i, \quad 1 \leq i,j \leq r\}.
     \endalign$$

We define $\tilde U_f$ to be the smallest reduced complex analytic subspace 
of $U \times \Bbb
P^{r-1}$ containing the graph $\Gamma (\phi_f)$.  
The support of $\tilde U_f$ is the closure of $\Gamma
(\phi_f)$ in the usual topology.  

The {\bf blow-up map} of $U$ along $\Cal I$ 
is the projection $\pi : \tilde U_f 
\rightarrow U$, which is a proper map.   

We will now show that the complex space 
$\tilde U_f$ is independent of the generators $f$ chosen for $I$. 

\proclaim {Lemma III.1}   
If $\{ f_1, ... , f_r \}$ and $\{ g_1, ..., g_s\}$ are two collections of generators of 
$\Cal I$ on $U$ then 
    $$\tilde U_f \cong \tilde U_g.$$
\endproclaim

\demo {Proof}   
Define a map $\psi : \Gamma(\phi_f) \rightarrow \Gamma(\phi_g)$ by
     $$\psi (z, [\xi]) = (z, [g_1(z): ... :g_s(z)].$$
The map $\psi$ is well-defined because $g_1(z), ... ,g_s(z)$ are not all $0$ for $z \in U -
V(I)$, (since $g_1, ... ,g_s$ are generators of $I$).  Furthermore $\psi^{-1}$ exists and is
given by
     $$\psi^{-1} (z,[\zeta]) = (z, [f_1 (z): ... :f_r(z)].$$
Both $\psi$ and $\psi^{-1}$ are clearly holomorphic so $\Gamma (\phi_f) \cong \Gamma
( \phi_g)$.  We will now show that they extend to holomorphic  maps on $\tilde U_f$ and
$\tilde U_g$. 

Since $\{ f_1, ... , f_r\}$ and $\{ g_1, ... , g_s \}$ both generate $I$, there exist
$\alpha_{ij}, \beta_{ij} \in \Cal O(U)$ such that
     $$g_i(z) 
     = \sum_{j=1}^r \alpha_{ij} (z) f_j (z) $$
and      $$f_i (z) 
     = \sum_{j = 1}^s \beta_{ij} (z) g_j (z) $$
for all $z$ in $U$.  Briefly, 
     $$f(z) = \beta (z) g(z) = \beta (z) \alpha (z) f(z) \tag *$$ 
for all $z \in U$.  The functions  
$\alpha$ and $\beta$ might not define maps on all of $\Bbb P^{r-1}$
and $\Bbb P^{s-1}$ but they do define maps on $\Gamma (\phi_f)$ and $\Gamma
(\phi_g)$.  

Suppose that $(z^\prime, [\xi^\prime]) \in U \times \Bbb P^{r-1}$ is the limit of 
points $(z_\gamma, [\xi_\gamma]) 
\in \Gamma (\phi_f)$, i.e. there is a sequence of points $\{
z_\gamma \} \in U$ such that 
     $$ z_\gamma \rightarrow z^\prime \qquad \text{and} \qquad 
   [\xi_\gamma ]= [f_1(z_\gamma) : ... :f_r (z_\gamma)] \rightarrow 
     [\xi^\prime].$$
Some component of $[\xi^\prime]$ is nonzero, say the first component, so that we may
assume that $\xi^\prime = (1, \xi_2^\prime , ... \xi_r^\prime)$.  Then we may also assume
that the sequence $\{ z_\gamma \}$ has the property that 
$f_1 (z_\gamma) \neq 0 $
for all $\gamma$ and that the sequence $\xi_\gamma$ is of the form 
     $$\xi_\gamma = (1, \xi_{\gamma 2}, ... ,\xi_{\gamma r})
     = \left( 1, \frac {f_2 (z_\gamma)}{f_1 (z_\gamma)}, ... , 
     \frac {f_r (z_\gamma)}{f_1 (z_\gamma)} \right) \tag {**}$$ 
where 
     $$\xi_\gamma  
       \rightarrow \xi^\prime.$$
We will use this description to show that $\alpha (z^\prime) \xi^\prime \neq 0$.   We
have 
     $$\alignat2 \beta (z_\gamma) \alpha (z_\gamma) \xi_\gamma 
     & = \beta (z_\gamma) \alpha (z_\gamma) \frac {f(z_\gamma)} {f_1 (z_\gamma)} 
     && \qquad \text{by (**)}  \\
     & = \frac {f(z_\gamma)} {f_1 (z_\gamma)} 
     && \qquad \text{by (*)} \\
     & = \xi_\gamma && \qquad \text{by (**).} \endalignat$$
Thus 
     $$\align \beta (z^\prime ) \alpha (z^\prime) \xi^\prime &= 
     \lim_{\gamma \rightarrow \infty}   \beta (z_\gamma) 
       \alpha (z_\gamma) \xi_\gamma \\
       &= \lim_{\gamma \rightarrow \infty} \xi_\gamma \\
       & = \xi^\prime \endalign$$
by continuity of $\alpha$ and $\beta$.
In particular, $\alpha (z^\prime) \xi^\prime \neq 0$ so $[\zeta] = [\alpha (z^\prime)
\xi^\prime]$ exists as a point of $\Bbb P^{s-1}$ (and 
is independent of the choices of representatives $\xi^\prime$ and $\xi_\gamma$). 
       
We define $\psi$ on $(z^\prime, [\xi^\prime])$ to be 
     $$\psi (z^\prime, [\xi^\prime]) = (z^\prime , [\alpha(z^\prime)  \xi^\prime]).$$
The definition of $\psi^{-1}$ is similar.  Clearly these extensions of $\psi$ and
$\psi^{-1}$ to the closures of $\Gamma (\phi_f)$ and $\Gamma (\phi_g)$ are
holomorphic and their compositions are the identity, so we obtain the required
isomorphism $\tilde U_f \cong \tilde U_g$.  
\qed
\enddemo   

\definition{Local Blow-up}  From the preceding lemma we see that it makes sense to
define the
blow-up of $U$ along
$\Cal I$ as $\text{Bl}_{\Cal I} U = 
\tilde U = \tilde U_f$ for any set of generators $f$.  
\enddefinition

If $\Cal I$ is the ideal of a {\bf smooth} 
subspace $C$ of $U$ then $\tilde U$ is also smooth.  
The set $C$ is called the {\bf center} of the blow-up.  
If $\Cal I$ is the ideal of a singular subset of
$U$ then $\tilde U$ may be singular.  

\proclaim{Lemma III.2}  Let $\Cal I$ and $\Cal J$ be nonzero 
coherent ideal sheaves on $U$ which are generated by global sections 
on $U$.  Suppose that $\Cal J$ is {\it principal}, i.e. 
generated by a single function on $U$.  Then  
     $$\text{Bl}_{\Cal I \Cal J} U \cong \text{Bl}_{\Cal I} U.$$
\endproclaim

\demo{Proof}  
Suppose that $\Cal J$ is generated locally by the single function $h$.  Then
     $$[hf_1: ... : hf_r] = [f_1: ... : f_r]$$ 
on $U - V(\Cal I \Cal J)$.  
\qed
\enddemo

We will use this lemma a little later to prove a corresponding statement about 
line bundles (Lemma III.4).  

\head Global Description of Blow-ups  \endhead

Let $\Cal I$ be a coherent sheaf of ideals on a complex manifold $M$.  
By Lemma III.1, we may extend the local definition of the blow-up canonically, to 
define a global blow-up 
     $$\pi : \tilde M = \text{Bl}_{ \Cal I} M 
       \rightarrow M.$$
The  blow-up map $\pi$  is proper and the 
restriction of $\pi$ from $\tilde M - \pi^{-1}(V(\Cal I))$ to $M - V(\Cal I)$ is biholomorphic.

If $\Cal I$ is
the ideal sheaf of a {\bf smooth} submanifold $C$ of $M$, then $\tilde M$ is smooth. 

\head Ideals, Divisors, Line Bundles, and Sections \endhead

Let $M$ be a complex manifold and let $D$ be a divisor on $M$.  
We denote by $L_D$ or  $[D]$ the corresponding line bundle on $M$.  
Let $\Cal L_D$ be the invertible sheaf of holomorphic sections of $[D]$.  

Let $s_D$ be a meromorphic section of $[D]$ whose divisor $(s_D)$ is $D$.  Such a
section always exists:  if $D$ is defined on an open covering $\{U_i\}$ of $M$ by
meromorphic functions $\{f_i\}$, the functions $\{f_i\}$ themselves define such a
section $s_D$. 

If $s$ is any other meromorphic section of $[D]$ then $ s \over {s_D}$ is a meromorphic
function on $M$.  Let $\Cal K_M$ be the sheaf of meromorphic functions on $M$.  We
may embed $\Cal L_D$ into $\Cal K_M$ by the map
     $$ s \mapsto {s \over {s_D}},$$
i.e. if $U$ is any open set in $M$ and $s \in \Cal L_D (U)$, we map $s$ to ${s \over
{s_D}} \in \Cal K_M (U)$. 

Now suppose that $Y$ is an effective divisor (codimension one
subscheme) of $M$ with ideal sheaf $\Cal I_Y$, and that $Y$ is given on an open cover
$\{ U_i \}$ of $M$ by holomorphic functions $\{ f_i\}$.  Let $s_Y$ be the
corresponding holomorphic section of $[Y]$.  Then $ 1 \over {s_Y}$ is a meromorphic
section of $[-Y]$.  We may embed $\Cal L_{-Y}$ into $\Cal K_M$ by the map
     $$s \mapsto  s s_Y.$$
The image of $\Cal L_{-Y}$ in $\Cal K_M$ is just the ideal $\Cal I_Y$ in $\Cal O_M
\subset \Cal K_M$.  Therefore
     $$\Cal L_{-Y} \cong \Cal I_Y.$$

Suppose that $\Cal I$ is any coherent ideal in $\Cal O_M$.  Tensoring the exact sequence
      $$\inj {\Cal I} {\Cal O_M}$$
by $\Cal L_{-Y}$ gives an exact sequence 
     $$\inj {\Cal I \otimes \Cal L_{-Y}} 
     {\Cal O_M \otimes \Cal L_{-Y} = \Cal L_{-Y}}$$ 
by Lemma II.3 above.  The image of $\Cal I \otimes \Cal L_{-Y}$ in $\Cal L_{-Y}$ is
just $\Cal I \Cal L_{-Y}$ (see e.g. [Ma],  p. 18). The image of $\Cal I \Cal L_{-Y}$ under
the embedding $\Cal L_{-Y} \hookrightarrow \Cal K_M$ is then $\Cal I \Cal I_Y$. 
Therefore

\proclaim{Lemma III.3}  Let $\Cal I$ be a coherent sheaf of ideals 
on a complex manifold $M$ and let $Y$ be an effective divisor on $M$.  Then 
    $$\Cal I \otimes \Cal L_{-Y} 
     \cong \Cal I \Cal I_Y.$$
\endproclaim

\proclaim{Lemma III.4}  
If $\Cal I$ is a coherent sheaf of ideals on a complex manifold $M$, 
and $\Cal L$  is the sheaf of holomorphic sections of a holomorphic line bundle on 
$M$,  then the blow-up of $M$ along $\Cal I$
is biholomorphic to the blow-up of $M$ along $\Cal I \otimes \Cal L$.   In particular, 
if $Y$ is an effective divisor on $M$, then the blow-up of $M$ along $\Cal I$ 
is biholomorphic to the blow-up of $M$ along $\Cal I \otimes \Cal L_{-Y} \cong 
\Cal I \Cal I_Y$.  
\endproclaim

\demo{Proof} 
Apply Lemma III.2. 
\qed
\enddemo

\proclaim{Lemma III.5}  
Let $M$ be a complex manifold and let $\Cal I$ be a coherent sheaf of ideals on $M$. 
Let $\pi : \tilde M = \text{Bl}_{ \Cal I} M \rightarrow M$ be the blow-up of $M$ along
$\Cal I$.  Then $\pi^{-1} \Cal I$ is a sheaf of principal ideals on $M$ (i.e. an invertible
sheaf).  The complex subspace of $\tilde M$ corresponding to $\pi^{-1} \Cal I$ is a
hypersurface. 
\endproclaim

\demo{Proof}  
Suppose that $\Cal I$ is generated locally on an open set $U$ in $M$ by $f_1, ... , f_r$.  Since
$\tilde U$ is contained in the subset of $U \times \Bbb P^{r-1}$ given by the equations
$f_i(z) \xi_j = f_j(z) \xi_i$, it is enough to prove that the inverse image ideal of $\Cal I$
on this set is principal.  But this is clear since on the set $U_i = \{ \xi_i \neq 0 \}$, we have 
     $$f_j = \frac {\xi_j} {\xi_i} f_i$$
so $f_i$ generates the inverse image ideal of $\Cal I$ on $U_i$. 
\qed
\enddemo

\head Exceptional Divisors of Blow-ups  \endhead

 The hypersurface in $\tilde M$ corresponding to $\pi^{-1} \Cal I$, 
described in Lemma III.5 above, is called the {\bf exceptional divisor} $E$ of $\pi$, i.e. 
     $$E = V(\pi^{-1} (\Cal I)) = \pi^{-1} V(\Cal I).$$

The proof of Lemma III.5 above gives us a local description of $E$.  
Suppose that $f_1, ... , f_r$ generate $\Cal I$ 
on an open set $U$ in $M$.   Cover $\tilde U \subset U \times
\Bbb P^{r-1}$ by sets $U_i = \{ \xi_i \neq 0 \}$. Then $E$ is given on $U_i$ by $f_i =
0$.

The map $\pi : \tilde M \rightarrow M$
is a proper map which is biholomorphic from $\tilde M - E$ to $M - V(\Cal I)$. If $\Cal
I$ is the ideal sheaf of a smooth center $C$, i.e. $\Cal I = \Cal I_C$, then $\tilde M$ is
smooth, $E = \pi^{-1} (C)$ is a smooth submanifold of $\tilde M$, and for each $p \in
C$ the inverse image $E_p = \pi^{-1} (p)$ is biholomorphic to $\Bbb P^{k-1}$, where
$k$ is the codimension of $C$ in $M$. 

\head Exceptional Line Bundles of Blow-ups  \endhead

Corresponding to the exceptional divisor $E$ on $\tilde M$ is an exceptional 
line bundle $L_E = [E]$.  Both $E$ and $L_E$ are independent of the 
local generators of $\Cal I$ used to construct the blow-up.  

In terms of local generators $f_1, ... , f_r$ of $\Cal I$, transition functions for $L_E$ are
     $$g_{ij} = { {f_i} \over {f_j}} = \frac {\xi_i} {\xi_j},$$
i.e. if $s$ is a holomorphic section of $L_E$ over $\tilde U$ then $s$ is represented by
holomorphic functions $s_i$ on $U_i = \{ \xi_i \neq 0 \}$ with
    $$s_i = g_{ij} s_j \qquad \text{on $U_i \cap U_j$.}$$

Since local transition functions for $L_E$ on the set $\tilde U$ 
are of the form $g_{ij} =  \frac {\xi_i} {\xi_j}$, 
the line bundle $L_E$ on $\tilde U$ is the restriction of the universal bundle $\Cal O (-1)$ 
on $U \times \Bbb P^{r-1}$.  More precisely, 
let $\sigma_1 : U \times \Bbb P^{r-1} \rightarrow U$ and $\sigma_2 : U \times \Bbb
P^{r-1} \rightarrow \Bbb P^{r-1}$ be the first and second projection maps, as shown below.   
     $$\CD
     \text{Bl}_I U = \tilde U
     @>>>
     U \times \Bbb P^{r-1} 
     @> \sigma_2 >>
     \Bbb P^{r-1} \\
     @.   @V  \sigma_1 VV   @.  \\
     @.
     U  @.
     \endCD $$
Let $\Cal O_{\Bbb P^{r-1}}(-1)$ be the universal bundle on $\Bbb
P^{r-1}$.  Then the restriction to $\tilde U$ of the line bundle 
$\sigma_2^* \Cal O_{\Bbb
P^{r-1}} (-1)$  is $L_E$ on $\tilde U$.  

We may interpret the fibre of $L_E$ over  $(z, [\xi]) \in \tilde U$ as the line through $\xi$ in
$\Bbb C^r$. 

\head Universal Property of Blow-ups \endhead

\proclaim{Lemma III.6 (Universal Property of Blow-ups)}  
Let $M$ be a complex manifold and let $\Cal I$ be a coherent sheaf of ideals on $M$. 
Let $\pi : \tilde M =  \text{Bl}_{ \Cal I} M \rightarrow M$ be the blow-up of $M$ along
$\Cal I$.  Suppose that $\phi : N \rightarrow M$ is a holomorphic map of a complex
space $N$ to $M$, such that the inverse image ideal $\phi^{-1} \Cal I$ is principal (i.e.
an invertible sheaf).  Then there exists a unique holomorphic lifting 
     $$\tilde \phi : N \rightarrow \tilde M$$
such that $\pi \circ \tilde \phi = \phi$.  
\endproclaim

\demo{Proof}  
Suppose that $f_1, ... ,f_r$ are generators
for $\Cal I$ over a small open set $U \subset M$.  Then $f_1 \circ \phi, ... , f_r \circ \phi$
are generators for $\phi^{-1} \Cal I$ over $\phi^{-1}(U)$ in $N$.  Since $\phi^{-1}
\Cal I$ is assumed to be 
a principal ideal sheaf,  all of the functions $f_i \circ \phi$
are multiples of one of them, so we have a well-defined map 
     $$\tilde \phi : \phi^{-1} (U) \rightarrow U \times \Bbb P^{r-1}$$ 
given by 
     $$v \mapsto (\phi(v) , [f_1 \circ \phi (v): ... : f_r \circ \phi (v)]).$$
By construction, the image of $\tilde \phi$ lies in the blow-up $\tilde U$ in $U \times
\Bbb P^{r-1}$ and $\pi \circ \tilde \phi (v) = \phi (v)$.  

By an argument 
similar to the proof of Lemma III.1 above, which showed that the blow-up $\tilde U$ is
independent of the collection of generators $\{ f_i \}$ used to construct it, we see that the
map $\tilde \phi$ is independent of the generators $\{ f_i \}$.
Thus we can extend our local construction to a well-defined holomorphic 
map $\tilde \phi : N \rightarrow \tilde M$.  

Finally we check the uniqueness of $\tilde \phi$.  Suppose that $\tilde \phi^\prime$ is any
holomorphic map from $N$ to $\tilde M$ such that 
$\pi \circ \tilde \phi^\prime = \phi =
\pi \circ \tilde \phi$.  
 Since $\pi$ is a biholomorphism away from the exceptional set, $\tilde
\phi^\prime$ and $\tilde \phi$ must agree on $\phi^{-1} (M - V( \Cal I)) = N - V(
\phi^{-1} \Cal I)$.  But $\phi^{-1} \Cal I$ was assumed to be a principal ideal, so
$V(\phi^{-1} \Cal I)$ is a hypersurface in $N$.  This means that $\tilde \phi^\prime$ and
$\tilde \phi$ agree on a dense set of $N$, so they must agree everywhere.
\qed
\enddemo

\head Blow-up of a Product of Ideals  \endhead

We will show that the blow-up of a product of two ideals looks like the composite of two 
blow-ups.  Since we have defined blow-ups
only for smooth manifolds, we will restrict ourselves to the case in which the blow-up 
along one ideal is smooth, for example  if that ideal is the ideal of a smooth
submanifold.  

\proclaim{Proposition III.7}  
Let $M$ be a complex manifold and $\Cal I_1$ and $\Cal I_2$ coherent sheaves of ideals
on $M$. 
Let  $\pi : \text{Bl}_{\Cal I_1} M \rightarrow M$ be the blow-up of $M$ 
along $\Cal I_1$ 
and suppose that the blow-up space $\text{Bl}_{\Cal I_1} M$ is smooth.   
Then 
     $$\text{Bl}_{\Cal I_1 \Cal I_2} M 
       \cong \text{Bl}_{\pi^{-1} \Cal I_2} \text{Bl}_{\Cal I_1}M,$$
i.e. the blow-up of $M$ along the product ideal $
\Cal I_1 \Cal I_2$ is isomorphic to the blow-up of $M$ along $\Cal I_1$ followed by the
blow-up along the inverse image ideal of $\Cal I_2$. 
\endproclaim

\demo{Proof}  
We will apply the universal mapping property of blow-ups (Lemma III.6).  
Let $N = \text{Bl}_{\pi^{-1} \Cal I_2} \text{Bl}_{\Cal I_1}M$ and 
let $\phi : N \rightarrow M$ be
the composite of the blow-up maps.  Then $\phi^{-1} \Cal I_1$ and $\phi^{-1} \Cal I_2$
are principal ideal sheaves on $N$ so $\phi^{-1} (\Cal I_1 \Cal I_2)$ is also principal. 
By the universal mapping  property, $\phi$ lifts to a holomorphic map $\tilde \phi : 
N \rightarrow \text{Bl}_{\Cal I_1 \Cal I_2} M $.    This map is a biholomorphism away from
the exceptional sets.  

Similarly, if $\psi : \text{Bl}_{\Cal I_1 \Cal I_2} M \rightarrow M$ is the blow-up 
of $M$ along $\Cal I_1 \Cal I_2$,
then $\psi^{-1} \Cal I_1$ is a principal ideal sheaf on $\text{Bl}_{\Cal
I_1 \Cal I_2} M$ and we can lift $\psi$ to a map 
$\psi_1 : \text{Bl}_{\Cal I_1 \Cal I_2} M \rightarrow \text{Bl}_{\Cal I_1} M$.
Next we check that $\psi_1^{-1} (\pi^{-1} \Cal I_2)$ is again a principal ideal
sheaf, so that we can lift $\psi_1$ to a map $\tilde \psi : \text{Bl}_{\Cal I_1 \Cal I_2} M
\rightarrow \text{Bl}_{\pi^{-1} \Cal I_2} \text{Bl}_{\Cal I_1}M = N.$  

Since the maps $\tilde \psi$ and $\tilde \phi$ are holomorphic everywhere and are
inverses of each other on open dense sets, they must be inverses 
of each other everywhere.  
\qed
\enddemo

\proclaim{Corollary III.8} 
Let $M$ be a complex manifold, $C$ a smooth center in $M$, 
and $\Cal I_C$ the ideal sheaf of $C$.  Then the blow-up of $M$ along $\Cal
I_C$ is isomorphic to the blow-up along $\Cal I_C^d$ for any integer $d > 1$, i.e. 
     $$\text{Bl}_{\Cal I_C} M \cong \text{Bl}_{\Cal I_C^d} M.$$
\endproclaim

\demo{Proof}  Apply Proposition III.7, noting that $\pi^{-1} \Cal I_C$ is 
principal and that blowing-up along a principal ideal sheaf leaves a space 
unchanged.  
\qed
\enddemo

\head Direct Images under Blow-up Maps  \endhead

We conclude section III by showing that the direct image of an ideal sheaf under 
a blow-up map is an ideal sheaf.  As always, we assume that the ideal sheaf $\Cal I$ 
for our blow-up is not the zero sheaf, so that $C = V (\Cal I)$ has 
codimension at least 1.  

\proclaim{Lemma III.9}  
Let $\pi : \tilde M \rightarrow M$ be the blow-up of a complex manifold $M$
along a coherent sheaf of ideals $\Cal I$ on $M$.  
Let $\Cal J$ be a sheaf of ideals on
$\tilde M$.  Then 
the direct image $\pi_* \Cal J$ is a sheaf of ideals on $M$.  
If $\Cal J$ is coherent then 
so is $\pi_* \Cal J$.  
\endproclaim

\demo{Proof}  
We wish to define a map $\pi_* \Cal J \rightarrow \Cal O_M$ and show that it is
injective.  To define a sheaf map $\pi_* \Cal J \rightarrow \Cal O_M$, it is enough to
define presheaf maps $\pi_* \Cal J (U) \rightarrow \Cal O_M (U)$ 
for all open sets $U$ in $M$. 
To show that a map of sheaves $\pi_* \Cal J \rightarrow \Cal O_M$  is
injective, it is enough to show that $\pi_* \Cal J(U) \rightarrow 
\Cal O_M (U)$ is injective for
all open sets $U$ in $M$. 

Recall that $\pi_* \Cal J(U) = \Cal J ( \tilde U)$, where 
$\tilde U = \pi^{-1} (U)$.  If $U$ does not intersect $C = V(\Cal I)$,  
then $\tilde U \cong U$ and $\pi_* \Cal J (U)$ may be
identified naturally as an ideal  in $\Cal O_M (U)$.  
Now suppose that $U$ does intersect $C$ and
consider $g \in \pi_* \Cal J (U) = \Cal J ( \tilde U )$. 
Let  $E $ be the 
exceptional divisor of $\pi$ in $\tilde M$.  
Since
     $$ \tilde U - \tilde U \cap E \cong U - U \cap C,$$
we may define a holomorphic function $G$ on $U - U \cap C$ whose pullback to
$\tilde U - \tilde U \cap E$ is $g$.  
For each $p \in U \cap C$, 
the fibre $\pi^{-1} (p)$ is compact, since $\pi$ is proper.  
Therefore $g$ is constant on $\pi^{-1} (p)$ 
and bounded on a neighborhood of $\pi^{-1} (p)$ in $\tilde U$.
Thus the function $G$ is locally 
bounded in $U$, so $G$ extends 
uniquely to a holomorphic function on $U$ by 
Riemann's Removable Singularity Theorem.  
Since $\pi^*G$ and $g$ are holomorphic on $\tilde U$ and 
equal on the dense set $\tilde U - \tilde U \cap E$, they must be equal on all 
of $\tilde U$, i.e. $\pi^*G = g$ on $\tilde U$.  For each
$g \in \Cal J (\tilde U)$ there is a unique such $G \in \Cal O_M (U)$, so we have a
well-defined map
     $$\pi_* \Cal J (U) \rightarrow \Cal O_M (U).$$
Clearly $G$ is identically zero if and only if $g$ is identically zero, so the map is
injective.

By the Direct Image Theorem, $\pi_* \Cal J$ is coherent if $\Cal J$ is, since $\pi$ 
is proper.   
\qed
\enddemo

\head IV.  Chow's Theorem for Ideals \endhead

This section is devoted to the proof of Chow's Theorem for Ideals, 
using the Direct Image Theorem. 
References for the usual Chow's theorem are [F] and [M].  
In section V we will state some applications to blow-ups. 

\proclaim{IV.1 Chow's Theorem for Ideals}  Let $U$ be an open neighborhood 
of $\{0\}$ in $\Bbb C^r$ and let  $X$ be an analytic subset of 
$U \times \Pn$.   Let $\Cal I$ be a coherent sheaf of ideals on
$X$.   Then $\Cal I$ is {\bf relatively algebraic} in the following sense:  $\Cal I$ is
generated (after shrinking $U$ if necessary) by a finite number of homogeneous
polynomials in homogeneous $\Pn$-coordinates, with analytic coefficients in
$U$-coordinates.   
\endproclaim

Since a sheaf on $X \subset U \times \Bbb P^n$ may be considered as a sheaf on 
$U \times \Bbb P^n$, we will ignore $X$ and prove the theorem for 
a coherent sheaf of ideals $\Cal I$ on $U \times \Bbb P^n$.  Although we have 
assumed that $U$ is an open neighborhood of $\{0\}$ in $\Bbb C^r$, the same 
methods could be used for any complex space $U$.  When we say that $\Cal I$ 
is generated by homogeneous polynomials in homogeneous $\Bbb P^n$-coordinates, 
we mean that the dehomogenizations of these polynomials generate the ideal 
locally.  We will show at the end of this section 
 that we may choose all the polynomial generators 
of $\Cal I$ to be of the same degree $d$, for $d$ sufficiently large. 

The usual Chow's theorem follows directly from Theorem IV.1:  
if $Y$ is an analytic subset of $U \times
\Pn$ and $\Cal I = \Cal I_Y$ is the ideal sheaf of $Y$ on $X = U \times \Bbb P^n$, 
then (after shrinking $U$ if necessary) $Y$ is cut out by 
a finite number of homogeneous polynomials 
in $\Pn$-coordinates with analytic coefficients in $U$-coordinates.  

\demo{Outline of Proof of Chow's Theorem for Ideals}  Let $\Ctn1$ be the blow-up of
$\Cn1$ at the origin and let $\sigma_1$ and $\sigma_2$ be the two projection maps of
$U \times \Ctn1$ as shown: 
     $$\CD   U \times \Ctn1 @> \sigma_2 >>
     U \times \Pn \\
     @V \sigma_1 VV @. \\
     U \times \Cn1 @. {}
     \endCD$$
The map $\sigma_2$ is flat since $U \times \Ctn1$ is a line bundle over $U \times \Pn$,
the product of the identity on $U$ with the universal line bundle on $\Pn$. Thus $\sinv 2
\Cal I = \sgst 2 \Cal I$ (Lemma II.4).  This inverse image ideal sheaf is coherent (see
facts on inverse image ideals, section II).  The map $\sigma_1$ is proper, so the direct
image $\Cal J = \sigma_{1*} (\sinv 2 \Cal I)$ is also coherent, by the Direct Image
Theorem.  Furthermore, $\Cal J$ is a sheaf of ideals on $U \times \Cn1$, not merely a
sheaf of modules, since $\sigma_1$ is a blow-up (Lemma III.9).   We will show 
(Lemmas IV.2 - IV.5) that 
$\Cal J$ is generated by homogeneous polynomials in $\Bbb C^{n+1}$-coordinates 
on a neighborhood of $(0,0)$, 
and that the corresponding polynomials in homogeneous $\Bbb P^n$-coordinates 
generate $\Cal I$.  

More specifically, let $x=(x_1, ... ,x_r)$ and $y=(y_0, ... ,y_n)$ be coordinates for $U$ and
$\Cn1$.  If $F(x, y)$ is a holomorphic
function in a neighborhood of $(0,0)$ in $U \times \Bbb C^{n+1}$ and $\lambda \in \Bbb C^*$, 
let $F^{(\lambda)}$ be the holomorphic function given by
     $$F^{(\lambda)} (x, y) = F(x, \lambda y).$$  
We first show (Lemma IV.2) that 
     $$F \in \Jx0 
     \Leftrightarrow F^{(\lambda)} \in \Jx0  
     \qquad \forall \lambda \in \Bbb C^*.$$

We use a corollary of Krull's Theorem to show that if $F^{(\lambda)} \in \Jx0$ for all
$\lambda \in \Bbb C^*$ then each homogeneous term in $y$ of $F(x, y)$ is in $\Jx0$
(Lemma IV.3).   

 It follows from Lemma IV.3 that $\Jx0$ is generated by a collection of 
homogeneous polynomials in $y$ with analytic coefficients in $x$.  
We then show that $\Jx0$ is generated by a finite number of these homogeneous
polynomials (Lemma IV.4). 

Finally we check that 
the same polynomials that generate $\Cal J = \sigma_{1*} (\sigma_2^{-1} \Cal I)$ 
over a neighborhood of $(0,0)$ in $U
\times \Cn1$, generate $\Cal I$ over a neighborhood of $\{0 \} \times \Pn \subset U
\times \Pn$ (Lemma IV.5). 
\qed
\enddemo

We will now prove Lemmas IV.2 - IV.5 to complete the proof of 
Chow's Theorem for Ideals.  
As above, let 
$x = (x_1, ... , x_r)$ and $y = (y_0, ... , y_n)$ be coordinates for $U \subset \Bbb
C^r$ and $\Cn1$, and let $F^{(\lambda)} (x, y) = F(x, \lambda y)$. 

\proclaim{Lemma IV.2}  
A holomorphic function $F$ is a section of $\Cal J = \sigma_{1*}( \sigma_2^{-1} \Cal I)$ 
on a neighborhood of $(0,0) \in U \times \Bbb C^{n+1}$ if and
only if $\Flam$ is a section of $\Cal J$ in a neighborhood of $(0,0)$ for each $\lambda
\in \Bbb C^*$. 
\endproclaim

\demo{Proof}  
A holomorphic function  is a section of $\Cal J = \sigma_{1*}( \sigma_2^{-1} \Cal I)$ on
a neighborhood of $(0,0)$ in $U \times \Bbb C^{n+1}$ if and only if 
its pullback by $\sigma_1$ is a section of $\sigma_2^{-1} \Cal I$ 
on a
neighborhood of $\sinv 1 (0,0) \cong \{ 0 \} \times \Pn$ 
in $U \times \tilde {\Bbb C}^{n+1}$.  
Suppose that $F$ is a section of $\Cal J$ on a neighborhood of $(0,0)$.  
To show that $\Flam$ is a
section of $\Cal J$ on a neighborhood of $(0,0)$, it is enough  to show
that $\sgst 1 \Flam$ is a section of $\sigma_2^{-1} \Cal I$ on a neighborhood of 
$p$ for each
$p \in \sinv 1 (0,0)$.  This reduces the proof to a simple calculation in local coordinates
near $p$ and $q = \sigma_2 (p)$. 

Choose homogeneous coordinates $[\xi_0: ... : \xi_n]$ on $\Pn$ such that 
the point $q = \sigma_2(p)$ in $U \times \Bbb P^n$ is given by $q  
= ( 0, [1: 0: ... : 0])$.  Let $W \subset \{ \xi_0 \neq 0 \} \subset \Pn $ be a
neighborhood of $[1: 0: ... :0]$  and let $w_i = { {\xi_i} \over {\xi_0}}$, for $1 \leq i
\leq n$, be nonhomogeneous coordinates for $W$.  The preimage $\sinv 2 (U \times W)
\cong U \times \Bbb C \times W$ is a neighborhood of $p$ 
in $U \times \tilde {\Bbb C}^{n+1}$ with coordinates $(x, y_0,
w) = (x_1, ... , x_r, y_0, w_1, ... , w_n)$ in which $p = (0,0,0)$.  
The maps $\sigma_1$ and $\sigma_2$ are given by 
     $$\sigma_1(x,y_0,w) = (x, y_0,y_0w) \qquad 
     \text{and} \qquad \sigma_2(x,y_0,w) = (x,w).$$

Since the ideal sheaf $\Cal I$ is coherent, $\Cal I$ is generated on a neighborhood of
$q$ by a finite collection of holomorphic functions $G_1, ... , G_s$.  The pullbacks $\sgst
2 G_1, ... , \sgst 2 G_s$ generate $\sigma_2^{-1} \Cal I$ on a neighborhood of $p$.  
Since $\sgst
1 F$ is a section of $\sigma_2^{-1} \Cal I$ on a neighborhood of $p$, there exist holomorphic
functions $A_1, ... , A_s$ on a neighborhood of $p$ such that
     $$\sigma_1^* F (x,y_0,w) 
     = \sum_{i=1}^s A_i (x,y_0,w) \sgst 2 G_i (x,y_0,w) .$$

Fix $\lambda \neq 0$.  Then for $y_0$ close enough to $0$, $(x,\lambda y_0,w)$ 
is in the domain of
the functions $\sgst 1 F$ and $A_1, ... , A_s$ and
     $$\align \sigma_1^* F^{(\lambda)} (x,y_0,w)
     &= \sgst 1 F (x,\lambda y_0,w) \\
     &= \sum_{i=1}^s A_i (x,\lambda y_0,w) \sgst 2 G_i (x,\lambda y_0,w) \\
     &= \sum_{i=1}^s A_i (x,\lambda y_0,w) G_i (x,w) \\
     &= \sum_{i=1}^s A_i (x,\lambda y_0,w) \sgst 2 G_i (x,y_0,w) . \endalign$$
Let $A_i^{(\lambda)}(x,y_0,w) = A_i (x,\lambda y_0,w)$ for $1 \leq i \leq s$.  Then each 
$A_i^{(\lambda)}$ is holomorphic on a
neighborhood of $p$ and
     $$\sgst 1 \Flam (x,y_0,w) 
     = \sum_{i=1}^s  A_i^{(\lambda)}(x,y_0,w) \sgst 2 G_i (x,y_0,w),$$ 
i.e. $\sgst 1 \Flam$ is a section of $\sigma_2^{-1}\Cal I$ on a neighborhood of $p$. 
\qed
\enddemo

\proclaim{Lemma IV.3} 
If $F^{(\lambda)} (x,y)$ is a section of $\Cal J$ on a neighborhood of $(0,0) 
\subset U \times \Bbb C^{n+1}$ for all
$\lambda \in \Bbb C^*$, then each homogeneous term in $y$ of $F(x, y)$ is a section of
$\Cal J$ on a neighborhood of $(0,0)$. 
\endproclaim

\demo{Proof}  
For any holomorphic function $F$ on a neighborhood of $(0,0)$,
let 
     $$F(x, y) = \sum_\alpha a_\alpha (x) y^\alpha$$
be the expansion of $F(x,y)$ in terms of monomials $y^\alpha = y_0^{\alpha_0}
y_1^{\alpha_1} ...  y_n^{\alpha_n}$ in $y$ with analytic coefficients $a_\alpha (x)$ in
$x$.  Let $\abs \alpha = \alpha_0 + \alpha_1 + ... + \alpha_n$.  The homogeneous term in
$y$ of degree $k$ in $F$ is
     $$F_k (x, y) = \sum_{\abs \alpha = k} a_\alpha (x) y^\alpha .$$ 
Then 
     $$F = \sum_{k=0}^\infty F_k 
     \qquad \text{and} \qquad 
      F^{(\lambda)} = \sum_{k=0}^\infty \lambda^k F_k .$$  
We wish to show that if $F$ is a section of $\Cal J$ on a neighborhood of $(0,0)$, then
each $F_k $ is also a section of $\Cal J$ on a neighborhood of $(0,0)$.  To minimize the
use of subscripts, we will also use $F$ and $F_k$ to represent the germs of these
functions at $(0,0)$. 

Let $A = \OXCx0$ (a Noetherian local ring), $(y) = (y_0, ... ,y_n)$ (an ideal contained in
the unique maximal ideal in $A$), and $J = \Jx0$ (also an ideal in $A$). Let
     $$\text{Jet}_m (F) = \sum_{k=0}^m F_k$$
be the $m$-jet of $F$ with respect to $y$.  Note that $F - \text{Jet}_m(F) \in
(y)^{m+1}$. 

By a corollary of Krull's Theorem (see e.g. [K],  Corollary 5.7,  p. 151), 
     $$J = \cap_{m \geq 0} (J + (y)^m),$$
where $(y)^0$ is defined to be $A$.  Since
     $$A = J +(y)^0 \supset J +(y)^1 \supset J+(y)^2 \supset ... $$
it follows that 
     $$J = \cap_{m \geq m_0} (J + (y)^m)$$
for any $m_0 \geq 0$.  

Suppose that $F^{(\lambda)} \in J$ for all $\lambda \in \Bbb C^*$.  Then since
     $$F^{(\lambda)} - \text{Jet}_m (F^{(\lambda)}) 
     \in (y)^{m+1}$$ 
we have  
     $$\text{Jet}_m (F^{(\lambda)}) \in J + (y)^{m+1}$$ 
for all $\lambda \in \Bbb C^*$.  Since $\text{Jet}_m (F^{(\lambda)}) = \sum_{k=0}^m
\lambda^k F_k$ for all $\lambda \in \Bbb C^*$, by taking $m + 1$ values of $\lambda$ 
it follows that
     $$F_k \in J + (y)^{m+1}$$ 
for $0 \leq k \leq m$.  
Fixing $k$, we have
     $$F_k \in J + (y)^{m+1} \qquad \text{for $m \geq k$}$$
or 
     $$F_k \in J + (y)^m \qquad \text{for $m \geq k+1$,}$$
i.e. 
     $$F_k \in \cap_{m \geq k + 1} (J + (y)^m) .$$ 
By the corollary of Krull's Lemma mentioned above, $F_k \in J$ for all $k$.
\qed  
\enddemo

\proclaim{Lemma IV.4}  
If $\Jx0$ is generated by a collection of elements of $\OXxpoly$ which are homogeneous
in $y$, then $\Jx0$ is generated by a {\bf finite} collection of elements of $\OXxpoly$ which
are homogeneous in $y$. 
\endproclaim 

\demo{Proof} 
Throughout the proof, whenever we refer to homogeneous
functions, we mean functions which are homogeneous in $y$.  
The ring $\OXCx0$ is Noetherian.    
As an ideal of $\OXCx0$, the ideal $\Jx0$ must be finitely generated,
but we want generators which are in $\OXxpoly$ and homogeneous.  In order to keep
track of the rings and ideals involved, we use the following notation: 
     $$\alignat2  A &= \OXCx0 
     &&\qquad \text{(a Noetherian ring)} \\
     B &= \OXxpoly 
     &&\qquad \text{(a Noetherian subring of $A$)} \\
     J & = \Jx0 
     &&\qquad \text{(an ideal in $A$)} \\
     J^\prime &= J \cap B 
     &&\qquad \text{(an ideal in $B$)}. 
     \endalignat$$

Suppose that there is a collection $H$ (perhaps infinite) of homogeneous generators of
$J$ over $A$ such that $H \subset B$.  Then $H \subset J^\prime$.  Since $J^\prime$ is
an ideal in $B$ and $B$ is Noetherian, there exists a finite set $H^\prime \subset B$ such
that $H^\prime $ generates $J^\prime$ over $B$. 

Each element of $H^\prime$ must be a linear combination of a finite number of
homogeneous generators in $H$.  Thus $J^\prime$ is generated over $B$ by a finite
number of homogeneous generators in $H$, i.e. we may choose $H^\prime$ to be a finite
set of homogeneous elements. 

Since $H \subset J^\prime$, $H^\prime$ also generates $H$ over $B$. Since $B \subset
A$, $H^\prime$ generates $H$ over $A$.  Finally, since $H^\prime$ generates $H$ over
$A$ and $H$ generates $J$ over $A$, $H^\prime$ generates $J$ over $A$, i.e. there
exists a finite set $H^\prime \subset \OXxpoly$ of homogeneous polynomials in $J$ such
that $H^\prime$ generates $\Jx0$ over $\OXCx0$. 
\qed
\enddemo

Note that each homogeneous element of $\OXxpoly$ is represented on a neighborhood
of $(0,0)$ by a homogeneous polynomial in $y$ with analytic coefficients in $x$. 

\proclaim{Lemma IV.5} 
The same polynomials that generate $\Cal J$ over a neighborhood of $(0,0)$ in $U
\times \Cn1$, generate $\Cal I$ over a neighborhood of $\{ 0 \} \times \Pn$ in $U \times
\Pn$. 
\endproclaim

\demo{Proof}
Suppose that $\Cal J$ is generated in a
neighborhood of $(0,0)$ by $F_1 (x,y), ... , F_s (x,y)$, where $F_i (x,y)$ is a
homogeneous polynomial of degree $d_i$ in $y$ with analytic coefficients in $x$.  We
will show that $\Cal I$ is generated on a neighborhood of $\{ 0 \} \times \Pn$ 
in $U \times \Bbb P^n$ by the
corresponding polynomials $F_i (x , \xi)$, where $[\xi] = [\xi_0: ... :\xi_n]$ are
homogeneous coordinates for $\Pn$.  More precisely, we will show that
$\Cal I$ is generated on a neighborhood of any point $q \in \{ 0 \} \times \Pn$ by
dehomogenizations of $F_1, ... , F_s$ near $q$.  

Choose homogeneous
coordinates $\xi$ on $\Pn$ such that $q = ( 0 , [1:0:...:0])$.  Nonhomogeneous
coordinates on the set $ W = \{ \xi_0 \neq 0 \} \subset \Bbb P^n$  are
$w_i = {{\xi_i} \over {\xi_0}}$ for $1 \leq i \leq n$.  We will check that $\Cal I$ is
generated in a neighborhood of $q$ by the polynomials
      $${ {F_i (x, \xi)} \over {\xi_0^{d_i}} } 
     = F_i \left( x, \frac {\xi} {\xi_0} \right)
     = F_i (x, 1, w_1, ... , w_n).$$

First we look at the maps $\sigma_1$ and $\sigma_2$ in local coordinates.  
We may use $(x,y_0,w)$ 
as local coordinates in $\sigma_2^{-1} (U \times W) \cong 
U \times \Bbb C \times W$.  Local coordinates for $U \times \Bbb
C^{n+1}$ are $(x,y_0,y_1, ... ,y_n)$, where $y_i = y_0 w_i$ for $1 \leq i \leq n$.  The maps
$\sigma_1$ and $\sigma_2$ 
are given by
$$\sigma_1(x,y_0,w) = (x, y_0,y_0w) \qquad 
     \text{and} \qquad \sigma_2(x,y_0,w) = (x,w).$$

Suppose that $G$ is a holomorphic section of $\Cal I$ on a neighborhood of $q$ in $U \times
\Bbb P^n$.  Then $\sigma_2^* G$ is a holomorphic section of $\sigma_2^{-1} \Cal I$ in a 
neighborhood of $\sigma_2^{-1} (q) = \{ (0,y_0,0) : y_0 \in \Bbb C \}$.  Since the homogeneous
polynomials $F_1, ... ,F_s$ generate $\Cal J = \sigma_{1*} (\sigma_2^{-1} \Cal I)$ on a
neighborhood of $(0,0) \in U \times \Bbb C^{n+1}$, their pullbacks $\sigma_1^* F_1, ... ,
\sigma_1^* F_s$ generate $\sigma_2^{-1} \Cal I$ on a neighborhood of $\sigma_1^{-1} 
(0,0) \in U \times \tilde {\Bbb C}^{n+1}$.  
In particular, there exist holomorphic functions  
$A_1, ... ,A_s$ on a neighborhood of 
the point $(x = 0, y_0=0,w=0)$ in $U \times \tilde \Bbb C^{n+1}$ 
such that 
     $$\sigma_2^* G(x,y_0,w) 
       = \sum_{i=1}^s A_i (x,y_0,w) \sigma_1^* F_i (x,y_0,w)$$
on that neighborhood.  But $\sigma_2^*G(x,y_0,w) = G(x,w)$ is
independent of the value of $y_0$ and $\sigma_1^* F_i(x,y_0,w) = F_i(x,y_0,y_0w) = y_0^{d_i}F_i(x,1,w)$
since $F_i$ is homogeneous of degree $d_i$ in $y$.  Therefore 
     $$ G(x,w) = \sum_{i=1}^s A_i(x,y_0,w)y_0^{d_i}F_i(x,1,w). $$
Choose some fixed nonzero value of $y_0$, close enough to $0$ 
that $(x,y_0,w)$ is in the domain of all the functions $A_i$ for 
$x$ and $w$ close enough to $0$.  Define  
     $$a_i(x,w) =  A_i(x,y_0,w) {y_0}^{d_i}.$$ 
Then 
     $$G(x,w) = \sum_{i=1}^s a_i(x,w) F_i(x,1,w).$$
Since the functions $a_i$ are holomorphic on a neighborhood of the point $q = (x=0,w=0)$, 
and the functions $F_i(x,1,w)$ are the local dehomogenizations of the homogeneous polynomials 
$F(x,\xi)$, we are done.  
\qed
\enddemo

This completes the proof of Chow's Theorem for Ideals.  We now show 
that the homogeneous polynomial generators of the ideal sheaf $\Cal I$ 
can be chosen to be of the same degree $d$, for large enough $d$.  

\proclaim {Corollary IV.6}  Let $U$ be an open neighborhood of $\{ 0 \}$ in 
$\Bbb C^r$ and let $X$ be an analytic subset of $U \times \Bbb P^n$.  
Let $\Cal I$ be a coherent sheaf of ideals on $X$.  
Then (possibly after shrinking $U$) 
there exists a positive integer $d_0$ such that for all 
$d \geq d_0$ the ideal $\Cal I$ is generated by a finite number of 
degree $d$ homogeneous polynomials in homogeneous $\Bbb P^n$-coordinates 
with analytic coefficients in $U$-coordinates.
\endproclaim

\demo{Proof}  As before, we may treat $\Cal I$ as a sheaf on $U \times 
\Bbb P^n$.  
By Chow's Theorem for Ideals, we may choose a finite collection 
of homogeneous polynomials generating $\Cal I$.  We wish to show that we 
can choose homogeneous polynomials which are all of the same degree.  
Suppose that $F_1, ... ,F_s$ are homogeneous polynomials of degrees $d_1, ... ,
d_s$ generating $\Cal I$ on $U\times \Bbb P^n$.  Let $d_0$ be any integer 
at least as large as the largest of $d_1, ... ,d_s$.  Then replace each 
$F_i$ with the set of all $\xi^\alpha F_i$ as $\xi^\alpha$ runs through all 
degree $d_0 - d_i$ monomials in homogeneous coordinates $[\xi] = [\xi_0: ... :\xi_n]$ 
on $\Bbb P^n$, i.e. use all monomials of the form $\xi_0^{\alpha_0} 
\xi_1^{\alpha_1} ... \xi_n^{\alpha_n}$ where $\alpha_0 + \alpha_1 + ... + \alpha_n 
= d_0 - d_i$. At every point in $U \times \Bbb P^n$, the dehomogenizations 
of the polynomials $\xi^\alpha F_i$ generate the same ideal as the dehomogenization 
of the polynomial $F_i$.  
\qed
\enddemo

Degree $d$ homogeneous polynomials on $\Bbb P^n$ may be viewed as 
sections of $\Cal O (d)$, the sheaf of holomorphic sections of 
the $d$th power of the hyperplane bundle on $\Bbb P^n$.  By abuse 
of notation, we will also use $\Cal O(d)$ to refer to the corresponding 
sheaf on $U \times \Bbb P^n$, obtained by pullback from $\Bbb P^n$ 
under the projection map $U \times \Bbb P^n \rightarrow \Bbb P^n$.  
If $\Cal I$ is a coherent sheaf of ideals on $U \times \Bbb P^n$, 
holomorphic sections of 
$\Cal I \otimes \Cal O (d)$ may be represented by homogeneous polynomials 
of degree $d$ in homogeneous $\Bbb P^n$-coordinates with 
analytic coefficients in $U$-coordinates, whose local 
dehomogenizations are sections of $\Cal I$.   

We can thus restate Corollary IV.6 as follows. 

\proclaim{Corollary IV.7}   Let $U$ be an open neighborhood of $\{ 0 \}$ in 
$\Bbb C^r$ and let $X$ be an analytic subset of $U \times \Bbb P^n$.  
Let $\Cal I$ be a coherent sheaf of ideals on $X$.  
Then (possibly after shrinking $U$) 
there exists a positive integer $d_0$ such that for all 
$d \geq d_0$ the ideal $\Cal I \otimes \Cal O(d)$ 
is generated by a finite number global sections on $X \subset 
U \times \Bbb P^n$.  
\endproclaim

\head V.   Chow's Theorem Applied to Blow-ups \endhead

In this section we consider 
consider some consequences of Chow's Theorem for Ideals 
for blow-ups.

\proclaim{Corollary V.1 (Blow-ups are Relatively Algebraic)}  
Let $\pi : \tilde M \rightarrow M$ be the blow-up of a complex manifold $M$ along 
a coherent sheaf of ideals $\Cal I$.  Then for each point $p$ in $M$, there exists 
a neighborhood $U$ of $p$ in $M$ and an embedding of $\tilde U = \pi^{-1} (U)$ 
into $U \times \Bbb P^{r-1}$, for some $r$, such that $\tilde U$ is cut out 
by a finite number of homogeneous polynomials in homogeneous $\Bbb P^{r-1}$-coordinates 
with analytic coefficients in $U$-coordinates.  Furthermore, we 
may choose all the homogeneous polynomial generators to be of the 
same degree $d$ for $d$ sufficiently large.  
\endproclaim

\demo{Proof}  Choose $U$ small enough that $\Cal I$ is generated by 
global sections $f_1, ... , f_r$ on $U$ and let $\tilde U \hookrightarrow 
U \times \Bbb P^{r-1}$ be the induced embedding.  
Then use 
Corollary IV.6 of Chow's Theorem for Ideals, with $X = U \times \Bbb P^{r-1}$ 
and the ideal $\Cal I = \Cal I_{\tilde U}$.
\qed
\enddemo

Now consider a coherent sheaf of ideals $\Cal J$ on $\tilde M$.  
Corollary IV.7 tells us that 
if $U$ is a small enough open set in $M$ and $d$ 
is a large enough positive integer,  
the sheaf $\Cal J \otimes \Cal O(d)$ is generated by a finite number 
of global sections 
on $\tilde U \subset U \times \Bbb P^{r-1}$.  
Recall from section III that the restriction of 
$\Cal O(d)$ to $\tilde U$ is just $\Cal L_{-E}^d$, the sheaf of 
holomorphic sections of the $d$th power of the dual 
of the exceptional line bundle.  From this observation 
and from Lemma III.3, we have 
     $$ \Cal J \otimes \Cal O(d) \cong
     \Cal J \otimes \Cal L_{-E}^d \cong
    \Cal J \Cal I_E^d.$$

\proclaim{Corollary V.2}  Let $\pi : \tilde M \rightarrow M$ be the blow-up of a complex
manifold $M$ along 
a coherent sheaf of ideals $\Cal I$ and let $E$ be the exceptional divisor of $\pi$.  Let $\Cal J$ 
be a coherent sheaf of ideals on $\tilde M$.  Then for each point $p$ in $M$ there exists 
a neighborhood $U$ of $p$ in $M$, an embedding of $\tilde U = \pi^{-1} (U)$ 
into $U \times \Bbb P^{r-1}$, for some $r$, and an integer $d_0$ such that 
the ideal $\Cal J \Cal I_E^d$ is generated 
by a finite number of global sections on $\tilde U$ for all $d \geq d_0$.  
\endproclaim

\demo{Proof}  Construct an embedding $\tilde U \hookrightarrow U \times 
\Bbb P^{r-1}$ using local generators of $\Cal I$, as usual.  
Then use Corollary IV.7 of Chow's Theorem for Ideals, 
with $X = \tilde U$ and the coherent sheaf of ideals 
$\Cal J \Cal I_E^d$ on $\tilde U$.  
\qed
\enddemo

Alternatively, the existence of these global generators 
over $\tilde U$ can
be proved using the positivity of the line bundle $L_E^{-1}$ along fibres of the map 
from $E$ to its image in $M$,  
as in Hironaka and Rossi [HR], using results of Grauert.  Except for
the use of the Direct Image Theorem, our method is more explicit.  We show not only
that global sections exist on $\tilde U$,  
but how they are related to homogeneous polynomials in $\Bbb
P^{r-1}$-coordinates generating $\Cal I$ locally. 

   In the special case of compact projective manifolds, these constructions can
be made global, using an ample line bundle on the original manifold. 

Applying the previous corollary and noting that homogeneous polynomials 
on $U \times \Bbb P^{r-1}$ determine hypersurfaces of $\tilde U$, we obtain 
the following.  

\proclaim{Corollary V.3}  Let $\pi : \tilde M \rightarrow M$ be the blow-up of a complex
manifold $M$ along 
a coherent sheaf of ideals $\Cal I$ and let $\Cal J$ 
be a coherent sheaf of ideals on $\tilde M$.  Then for each point $p$ in $M$ there exists 
a neighborhood $U$ of $p$ in $M$, such that the complex space $V(\Cal J)$ determined 
by $\Cal J$ is cut out by a finite number of hypersurfaces in $\tilde U = \pi^{-1} (U)$. 
In particular, if $C$ is a smooth center in $\tilde M$ and $\Cal J = \Cal I_C$, 
then $C$ is cut out by hypersurfaces, not only locally in $\tilde M$, but 
over the pre-images $\tilde U$ of small open sets $U$ in $M$.  
\endproclaim

The next corollary will be instrumental in constructing single-step 
blow-ups.  

\proclaim{Corollary V.4}  Let $\pi : \tilde M \rightarrow M$ be the blow-up of a {\bf compact} 
complex manifold $M$ along 
a coherent sheaf of ideals $\Cal I$ and let $E$ be the exceptional divisor of $\pi$.  Let $\Cal J$ 
be a coherent sheaf of ideals on $\tilde M$.  Then there exists an integer $d_0$ such that 
     $$\pi^{-1} \pi_* ( \Cal J \Cal I_E^d) 
     = \Cal J \Cal I_E^d$$
for all $d \geq d_0$.  
\endproclaim

\demo{Proof} 
By compactness it is enough to prove the statement locally over neighborhoods of points in 
$M$.  
By Corollary V.2, for each point $p$ in $M$ there exists a neighborhood $U$,  
an embedding $\tilde U \hookrightarrow U \times \tilde \Bbb P^{r-1}$, for some $r$,  and an 
integer $d_0$ such that $\Cal J \Cal I_E^d$ is generated by a finite number of global 
sections on $\tilde U$, for $d \geq d_0$.  These sections are holomorphic functions, 
vanishing on $E$ for $d > 0$.  By the Riemann Extension Theorem, they determine holomorphic 
functions on $U$.  These functions on $U$ generate $\pi_* ( \Cal J \Cal I_E^d)$ 
and their pullbacks to $\tilde U$ generate $\pi^{-1} \pi_* ( \Cal J \Cal I_E^d)$.
Therefore $\pi^{-1} \pi_* (\Cal J \Cal I_E^d) = 
\Cal J \Cal I_E^d$.  
\qed
\enddemo

\example{Remark V.5}  Using local coordinates and local generators of 
$\Cal I$, we can describe more 
concretely the relationship between homogeneous polynomials 
generating $\Cal J $ over $\tilde U$ and 
holomorphic functions generating $\Cal J \Cal I_E^d$ over $\tilde U$.  

Since $\Cal I$ is coherent, $\Cal I$ is generated by a finite 
collection of holomorphic functions $f_1, ... , f_r$ on $U$, 
for $U$ small enough. 
Let $z$ represent $U$-coordinates 
and $[\xi] = [\xi_1: ... :\xi_r]$ homogeneous $\Bbb P^{r-1}$-coordinates.  
By Chow's Theorem for Ideals, $\Cal J$ is generated by a finite 
collection of homogeneous polynomials $F(z,\xi)$ (homogeneous in $\xi$ 
and analytic in $z$).  
 The ideal sheaf $\Cal I_E$ of the 
exceptional divisor is generated by the 
pullbacks of $f_1, ... ,f_r$ to 
$\tilde U$.  For simplicity we will also refer to these 
pullbacks as $f_1, ... , f_r$.  
The sheaf $\Cal I_E^d$ is generated by 
all monomials of degree $d$ in $f_1, ... , f_r$.  
The sheaf $\Cal J \Cal I_E^d$ is generated by all products of 
the form $f^\alpha F(z,\xi) $, where $f^\alpha$ 
represents a degree $d$ monomial in $f_1, ... , f_r$.   
The function $F(z,\xi)$ is of the form
     $$F(z,\xi) = \sum_\beta c_\beta (z) \xi^\beta$$
where $\xi^\beta$ is a monomial of degree $d$ in $\xi_1, ... ,\xi_r$ 
and $c_\beta (z)$ is a holomorphic function of $z$.  
Then 
     $$\alignat2 f^\alpha F(z, \xi) 
       &= \sum_\beta c_\beta (z) \xi^\beta f^\alpha &&\\
       & = \sum_\beta c_\beta (z) \xi^\alpha f^\beta 
       &&\qquad \text{since $f_i \xi_j = f_j \xi_i$.}
       \endalignat$$
Thus 
     $$ f^\alpha F(z, \xi) = \xi^\alpha F(z, f).$$
The sheaf $\Cal J \Cal I_E^d$ is generated by all such products as 
$\xi^\alpha$ ranges over all degree $d$ monomials in $\xi_1, ... 
, \xi_r$.  Since these monomials in $\xi$ cannot all be zero 
simultaneously, the collection $\{ \xi^\alpha F(z, f) 
\}_\alpha$ is generated by $F(z,f)$.  

     We now see explicitly the holomorphic generators of 
$\Cal I \Cal I_E^d$ described in the previous corollary - they are 
the functions $F(z,f)$.  These functions are holomorphic on  
$\tilde U$ and  
vanish on $E$ for $d > 0$, so they define holomorphic functions on $U$.  
As functions 
on $U$, they generate $\pi_* ( \Cal J \Cal I_E^d)$.  Their pullbacks to 
$\tilde U$ generate $\pi^{-1} \pi_* ( \Cal J \Cal I_E^d)$ and 
are once again the functions $F(z, f)$.
\endexample

\example{Example V.6}  
Let $\Cal I$ be ideal sheaf of the origin in $\Bbb C^3$
(i.e. $V (\Cal I) =
C = \{Z_1 = Z_2 = Z_3 = 0\}$), let $\pi : {\tilde {
\Bbb C}}^3 \rightarrow \Bbb C^3$ be the blow-up along $\Cal I$, and let $E = \pi^{-1} (C)$ 
be the exceptional divisor. 
Let $\Cal J$ be the ideal on ${\tilde {\Bbb C}}^3$ generated by the homogeneous
polynomial $\xi_1 \xi_2 - \xi_3^2$.  Let $F(Z) = Z_1 Z_2 - Z_3^2$ be the corresponding
polynomial on $\Bbb C^3$.  Then $\pi^* F$ is a holomorphic section of $\Cal J \Cal
I_E^2$. We have
     $$\Cal J \supset \Cal J \Cal I_E 
     \supset \Cal J \Cal I_E^2 \supset ...$$ 
and 
     $$\pi^{-1} \pi_* (\Cal J \Cal I_E^d) 
     = \cases \Cal J \Cal I_E^2 \qquad d < 2 \\
     \Cal J \Cal I_E^d \qquad d \geq 2. \endcases$$
Note that although we refer to $\xi_1 \xi_2 - \xi_3^2$ as a generator of $\Cal J$, it is not
a function on $ {\tilde {\Bbb C}}^3$.  If $U$ is any neighborhood of $0$ in $\Bbb
C^3$, the only nonzero holomorphic sections of $\Cal J$ on $\tilde U =
\pi^{-1} (U)$ are those
generated by homogeneous polynomials of degree at least 2, which must be vanishing on
$E$ to degree at least 2. 
\endexample

Once again, 
the next result could be proved using properties of positive line bundles 
with methods similar to those of Hironaka and Rossi in [GR] and results 
of Grauert.    In the algebraic setting it could be proved using ample line
bundles.   We restrict ourselves to the case in which the blow-up $\tilde M$ is smooth,  
since this is the only case we require and since we have defined the blow-up of $\tilde M$ along 
$\Cal J$ only in the case in which $\tilde M$ is smooth.  

\proclaim{Corollary V.7} 
Let $\pi : \tilde M \rightarrow M$ be the blow-up of a compact complex 
manifold along a coherent sheaf of ideals $\Cal I$ such that $\tilde M$ is smooth, 
and let $E$ be the exceptional divisor of $\pi$.    
Let $\Cal J$ be a coherent sheaf of ideals on $\tilde M$.  Then there exists an 
integer $d_0$ such that 
the blow-up 
of $\tilde M$ along $\Cal J$ 
is isomorphic to  the blow-up of $\tilde M$ along $\pi^{-1} \pi_* ( \Cal J \Cal
I_E^d)$ for all $d \geq d_0$.   
\endproclaim

\demo{Proof}  By Corollary V.4 there exists a $d_0$ such that 
 $\pi^{-1} \pi_* ( \Cal J \Cal I_E^d) = \Cal J \Cal I_E^d$ for all $d \geq d_0$.  
By Lemmas III.3 and III.4,  
the blow-up along $\Cal J$ is isomorphic to the blow-up along $\Cal J \Cal I_E^d$.
\qed
\enddemo

The direct image of a product is not always the product of the direct images.  
In the next lemma we give a condition under which products of ideal sheaves 
behave well under direct images of blow-up maps.  

\proclaim{Lemma V.8} 
Let $\pi : \tilde M \rightarrow M$ be the blow-up of a 
compact complex manifold $M$ along a coherent sheaf of ideals $\Cal I$ 
and let $E$ be the exceptional divisor.  
Let $\Cal J_1$ and $\Cal J_2$ be coherent sheaves of ideals 
on $\tilde M$.  Then for $d_1$ and $d_2$ large enough, 
     $$\pi_* (\Cal J_1 \Cal J_2 \Cal I_E^{d_1 + d_2}) = 
       \pi_* (\Cal J_1 \Cal I_E^{d_1})        \pi_* (\Cal J_2 \Cal I_E^{d_2}).$$
\endproclaim

\demo{Proof}  Since $M$ is compact, 
it is enough to prove the lemma locally, on a blow-up 
$\pi : \tilde U \rightarrow U$ of an open set $U$.  
We use the notation of remark V.5 above.  
By Corollary IV.6, if $\Cal J$ is 
a coherent sheaf of ideals on $\tilde U$, then for $d$ large enough 
and possibly after shrinking $U$, 
the ideal $\Cal J$ is generated on $\tilde U \subset U \times \Bbb P^{r-1}$ 
by a finite number of degree $d$ homogeneous 
polynomials $F(z,\xi)$ in homogeneous coordinates $\xi$ on $\Bbb P^{r-1}$.  
As was shown in remark V.5, the functions $F(z, f)$ 
generate the direct image $\pi_* ( \Cal J \Cal I_E^d)$.  

     If a finite collection $\{ F(z, \xi)\}$ of degree $d_1$ polynomials 
generates 
$\Cal J_1$ and a finite collection $\{ G(z, \xi) \}$ of degree $d_2$ polynomials 
generates $\Cal J_2$, then the collection $\{ F(z, f ) \}$ generates 
$\pi_* ( \Cal J_1 \Cal I_E^{d_1})$ and the collection $\{ G(z, f) \}$ 
generates $\pi_* (\Cal J_2 \Cal I_E^{d_2})$.  The collection of 
all products $F(z, f) G(z, f)$ generates $\pi_* ( \Cal J_1 
\Cal I_E^{d_1}) \pi_* (\Cal J_2 \Cal I_E^{d_2})$.  Similarly, 
the collection of all products $F(z, \xi) G(z, \xi)$ generates 
$\Cal J_1 \Cal J_2$, and since these products are degree $d_1 + d_2$ 
homogeneous polynomials in $\xi$, the collection of all products 
$F(z, f) G(z, f)$ generates $\pi_* ( \Cal J_1 \Cal J_2 
\Cal I_E^{d_1 + d_2})$.  Thus 
$\pi_* (\Cal J_1 \Cal J_2 \Cal I_E^{d_1 + d_2}) = 
\pi_* (\Cal J_1 \Cal I_E^{d_1})        \pi_* (\Cal J_2 \Cal I_E^{d_2})$.
\qed
\enddemo

\example{Remark V.9}
To see that the direct image of a product is not always the product of the 
direct images,  we refer to
Example V.6.  In that example, we described a sheaf of ideals $\Cal J$ 
on $\tilde {\Bbb C}^3$ generated by a degree $2$ homogeneous polynomial and such that 
     $$\pi^{-1} \pi_* ( \Cal J \Cal I_E^d) = 
     \cases  
     \Cal J \Cal I_E^2 & \quad d < 2 \\
     \Cal J \Cal I_E^d & \quad d \geq 2. \endcases$$

Suppose that $\pi_* (\Cal J \Cal I_E^d) = (\pi_* \Cal J)( \pi_* \Cal I_E^d)$.  Then 
     $$\alignat2 \pi^{-1} \pi_* (\Cal J \Cal I_E) 
       &= (\pi^{-1} \pi_* \Cal J)( \pi^{-1} \pi_* \Cal I_E ) 
     &&\qquad \text{by Lemma II.10} \\
     & = (\Cal J \Cal I_E^2) \Cal I_E &&\\ 
     & = \Cal J \Cal I_E^3 &&\endalignat$$
which is impossible since 
     $$\pi^{-1} \pi_* ( \Cal J \Cal I_E) = \ \Cal J \Cal I_E^2$$
by the example.
\endexample

\head VI.  Replacing a Sequence of Blow-ups by a Single Blow-up \endhead

Let $X$ be a singular subvariety of a compact complex manifold $M$.  
In this section we show how to replace a sequence of blow-ups along 
smooth centers, which resolves the singularities of $X$, by a single 
blow-up of $M$ along a coherent sheaf of ideals $\Cal I$, 
which is a product of coherent ideals corresponding to the centers.  
The support of $\Cal I$ is the singular locus of $X$, the proper 
transform of $X$ in the blow-up of $M$ along $\Cal I$ is nonsingular, 
and the exceptional divisor of the blow-up along $\Cal I$ is a 
normal crossings divisor which has normal crossings with the desingularization 
of $X$.  

\proclaim{Proposition VI.1}  Let $M$ be a compact complex manifold and let 
     $$M^{\prime \prime} \overset {\pi^\prime} \to \rightarrow
     M^\prime \overset \pi \to \rightarrow M$$
be a sequence of blow-ups such that 
 \roster
   \item"{a.}" $\pi : M^\prime \rightarrow M$ is the blow-up 
of $M$ along a coherent sheaf of ideals $\Cal I$ such that 
$M^\prime$ is smooth and 
$V(\Cal I)$ has codimension at least 2 and 
   \item"{b.}"  
     $\pi^\prime : M^{\prime \prime} \rightarrow M^\prime$     
is the blow-up of $M^\prime$ along a smooth center $C$ of codimension 
at least $2$.
 \endroster
Let $E$ be the exceptional divisor of $\pi$ in $M^\prime$.
Then the sequence of blow-ups $M^{\prime \prime} \rightarrow M^\prime 
\rightarrow M$ is equivalent to a single blow-up along a coherent 
sheaf of ideals $\Cal J$ on $M$ given by 
     $$\Cal J = \Cal I \Cal I^\prime$$
where $\Cal I^\prime = \pi_* (\Cal I_C \Cal I_E^d)$ and $d$ 
is a large enough positive integer that $\pi^{-1} \pi_* ( \Cal I_C 
\Cal I_E^d) = \Cal I_C \Cal I_E^d$.  Furthermore 
\roster
   \item"{i.}"  the blow-up of $M^\prime$ along $\pi^{-1} \Cal I^\prime 
= \Cal I_C \Cal I_E^d$ is isomorphic to the blow-up along $C$, i.e. 
the blow-up of $M^\prime$ along $\pi^{-1} \Cal I^\prime$ is 
isomorphic to $M^{\prime \prime}$, and 
     \item"{ii.}"  the complex space $V(\Cal J)$ determined by $\Cal J$ 
has codimension at least 2 in $M$.  
\endroster  
\endproclaim

\demo{Proof}  By Corollary V.4 
     $$\pi^{-1} \pi_* (\Cal I_C \Cal I_E^d) = \Cal I_C \Cal I_E^d$$
for all sufficiently large $d$.  We apply Proposition III.7 to 
$\Cal J = \Cal I \Cal I^\prime = \Cal I \ \pi_* ( \Cal I_C 
\Cal I_E^d)$ to show that blowing up $M$ along $\Cal J$ is equivalent 
to first blowing up $M$ along $\Cal I$ to obtain $M^\prime$, and then blowing 
up $M^\prime$ along $\Cal I_C \Cal I_E^d$.  But the blow-up  
along $\Cal I_C \Cal I_E^d$ is equivalent to the blow-up along 
$\Cal I_C$ by Lemma III.4.  

    Finally we note that 
     $$\align V(\Cal J) & = V(\Cal I) \cup V(\pi_* (\Cal I_C \Cal I_E^d))  
      \\
     &= V(\Cal I) \cup \pi ( V(\Cal I_C) \cup V(\Cal I_E^d))   \\
     &= V(\Cal I) \cup \pi(C)  \endalign$$
which has codimension at least 2.
\qed
\enddemo

We apply Proposition VI.1 inductively to obtain 

\proclaim{Proposition VI.2}  Let $M_0$ be a compact complex manifold and let 
     $$M_m \overset {\pi_m} \to \rightarrow M_{m-1} 
      \overset {\pi_{m-1}} \to \rightarrow ... 
     \overset {\pi_2} \to \rightarrow M_1  
     \overset {\pi_1} \to \rightarrow M_0$$
be a sequence of blow-ups along smooth centers $C_j \subset M_{j-1}$ 
of codimension at least 2.   Then the composite $\pi_1 \circ ... \circ \pi_m 
: M_m \rightarrow M_0$ is equivalent to a single blow-up 
 along a coherent sheaf of ideals 
     $$\Cal I = \Cal I_1 \Cal I_2 ... \Cal I_m$$
where $\Cal I_1, \Cal I_2, ... , \Cal I_m$ are coherent sheaves 
of ideals on $M$ such that  
\roster 
   \item"{i.}"  the blow-up of $M_{j-1}$ along the inverse image 
ideal of $\Cal I_j$ on $M_{j-1}$ is isomorphic to the blow-up of 
$M_{j-1}$ along $ C_j$, and  
   \item"{ii.}"  the complex space $V(\Cal I)$ has codimension 
at least 2 in $M_0$.  
\endroster
\endproclaim

\demo{Proof}  We construct the ideal sheaves $\Cal I_1, ... , \Cal I_m$ 
inductively, using Proposition VI.1, and noting that all the 
spaces $M_j$ are smooth, since the centers of the blow-ups are 
smooth.  We may construct an ideal sheaf $\Cal I_j$ 
from $\Cal I_{C_j}$ either step-by-step, going down one level at 
a time, or all in one step, using the composite of the first 
$j-1$ blow-ups.    We use  
the second method in this proof, because it is notationally simpler.  
The first method is computationally simpler, so we use it in 
our example in section IX.  

Start by letting 
$\Cal I_1  = \Cal I_{C_1}$, the ideal sheaf of the first center $C_1$, 
and construct $\Cal I_2$ as in Proposition VI.1.   The blow-up of 
$M_1$ along $\pi_1^{-1} \Cal I_2$ is isomorphic to $M_2$ and the 
complex space $V(\Cal I_1 \Cal I_2)$ has codimension at least 2.  
Next suppose that we have  constructed $\Cal I_1, ... , \Cal I_{j-1}$ 
satisfying condition (i), and such that $V(\Cal I_1 ... \Cal I_{j-1})$ 
has codimension at least 2 in $M_0$.  Condition (i) implies that 
the blow-up of $M_0$ along  the product $\Cal I_1 ... \Cal I_{j-1}$ 
is isomorphic to $M_{j-1}$.  Let 
     $$\tau = \pi_1 \circ ... \circ \pi_{j-1} : M_{j-1} \rightarrow M_0$$ 
be this blow-up map and let $D$ be the exceptional divisor of $\tau$ 
in $M_{j-1}$.  Pick $d$ large enough such that $\tau^{-1} \tau_* 
( \Cal I_j \Cal I_D^d) = \Cal I_j \Cal I_D^d$  
and set    
     $$\Cal I_j = \tau_* ( \Cal I_j \Cal I_D^d).$$
Then apply Proposition VI.1.
\qed
\enddemo

Using Hironaka's theorem on the existence of embedded resolutions of 
singularities we obtain 

\proclaim{Corollary VI.3}  Let $M$ be a compact complex manifold and 
let $X$ be a singular subvariety of $M$.  Let 
$$M_m \overset {\pi_m} \to \rightarrow M_{m-1} 
      \overset {\pi_{m-1}} \to \rightarrow ... 
     \overset {\pi_2} \to \rightarrow M_1  
     \overset {\pi_1} \to \rightarrow M_0=M$$
be a sequence of blow-ups along smooth centers $C_j \subset M_{j-1}$ 
of codimension at least 2 which resolves the singularities of 
$X$, and such that the total exceptional divisor of the composite 
map has normal crossings and has normal crossings with the desingularization 
of $X$ in $M_m$.  
Then there exists a coherent sheaf of ideals 
$\Cal I$ on $M$ of the form 
     $$\Cal I = \Cal I_1 \Cal I_2 ... \Cal I_m$$
such that for each $j$, the blow-up map of $M$ along $\Cal I_1 \Cal I_2 ... \Cal I_j$ 
is equivalent to the composite map $\pi_1 \circ \pi_2 \circ ... \circ 
\pi_j : M_j \rightarrow M_0$.  In particular, 
 \roster
   \item"{i.}"  the proper transform $\tilde X$ of $X$ in the blow-up $\tilde M$ 
of $M$ along $\Cal I$ is nonsingular,
   \item"{ii.}"  $V(\pi^{-1} \Cal I)$ is a normal crossings divisor in $\tilde M$ 
which has normal crossings with $\tilde X$, and 
   \item"{iii.}"  the support of $\Cal I$ is the singular locus of $X$ in $M$.
\endroster
\endproclaim

\head VII. Chern Forms and Metrics for Exceptional Line Bundles  \endhead

Let $\pi : \tilde M \rightarrow M$ be the blow-up of a compact complex 
manifold $M$ along a coherent sheaf of ideals $\Cal I$ 
such that $\tilde M$ is smooth.  Let $E$ be 
the exceptional divisor of $\pi$ and $L_E = [E]$ the associated line 
bundle on $\tilde M$.  In this section we describe explicitly the construction 
of a Chern form on $L_E$ which is 
negative definite on the fibres of the map $E \rightarrow C = V(\Cal I)$.  

We first construct local Chern forms on  sets of the form $\tilde U = \pi^{-1} 
(U)$, where $U$ is a small open set in $M$.  An embedding 
$\tilde U \hookrightarrow U \times \Bbb P^{r-1}$ induces a local metric and 
local Chern form on the line bundle $L_E$ over $\tilde U$,  using
the Fubini-Study form on $\Bbb P^{r-1}$.  
Different embeddings of $\tilde U$
corresponding to different choices of local generators of $\Cal I$ may give different
Chern forms in the same Chern class.
This type of local Chern form has a particularly simple formula 
in terms of the local generators of $\Cal I$.  It 
is negative definite on the fibres of the 
map  $E \rightarrow C$ and negative semi-definite on 
$\tilde U$, since it is the pullback of the
negative of the Fubini-Study form on $\Bbb P^{r-1}$. 
 We then patch globally
using $C^\infty$ partitions of unity on $M$, 
to obtain  global metrics and Chern forms for $L_E$. 

\head Chern Forms on Line Bundles \endhead

We begin with some background material on Chern forms. 
Let  $L \rightarrow N$ be a holomorphic line bundle on a complex 
manifold $N$.  Choose a cover of $N$ by open sets $V_i$ such that 
$L$ is trivial on $V_i$,  
and let $\{ g_{i j }\}$ be holomorphic transition 
functions for a trivialization of $L$ over $\{ V_i \}$.  
A holomorphic section $s$ of $L$ over 
$N$ may be given by a collection of 
holomorphic functions $s_i$ on $V_i$ which transform on $V_i 
\cap V_j$ by the rule       
   $$s_i = g_{i j} s_j.$$
A hermitian metric $h$ on $L$ may 
be described by a collection of positive $C^\infty$ functions $h_i$ on 
$V_i$ such that the norm of $s$ is given on $V_i$ by 
     $$\normsq s = \absq {s_i} h_i.$$
The functions $h_i$ transform by the rule 
     $$h_j = \absq {g_{i j}} h_i.$$

\subhead  Local description of a Chern form \endsubhead 
The Chern form of $L$ with respect to $h$ is given on $V_i$ by 
     $$c_1(L,h) = - \rn2pi \ppbar \log h_i.$$
Note that this (1,1)-form is well-defined on $N$, because  
     $$\alignat2 \ppbar \log h_j 
       & = \ppbar \log \absq {g_{i j }} h_i &&\\
       &= \ppbar ( \log g_{i j} + \log {\overline g}_{i j} 
       + \log h_i )&&\\
       &= \ppbar \log h_i 
       &&\qquad \text{since $g_{i j}$ is holomorphic.} \endalignat$$

\subhead  Formula for a Chern form off the zero locus of a section $s$ \endsubhead 
On the set on which $s \neq 0$ we may write 
     $$c_1(L,h) = - \rn2pi \ppbar \log \normsq s.$$
                                 
\subhead  Chern forms via pullbacks  \endsubhead  
Chern forms behave well under pullbacks.  Suppose that $\phi : N_1 \rightarrow N_2$ 
is a holomorphic map of complex manifolds and $L$ is a line bundle on $N_2$ 
with metric $h$. Then $\phi^* L$ is a line bundle on $N_1$ with an induced metric 
$\phi^* h$, and the Chern form of $\phi^* L$ with respect to $\phi^* h$ is
the pullback of the Chern form of $L$ with respect to $h$, i.e. 
     $$c_1 ( \phi^*L, \phi^* h) = \phi^* c_1 ( L,h).$$

\head   Local Chern Forms for Blow-ups \endhead

Let $U$ be 
an open set in $\Bbb C^n$ and let 
$\pi : \tilde U \rightarrow U$ be the blow-up of $U$ 
along a coherent sheaf of ideals $\Cal I$ such that 
$\tilde U$ is smooth.  We will assume that $U$ is small 
enough that $\Cal I$ is generated by global sections 
$f_1, ... , f_r$ over $U$.  Let $E$ be the exceptional divisor 
and $L_E$ the associated line bundle on $\tilde U$.  

If $\Cal I$ is generated by a single function over 
$U$, then the sheaf $\Cal I$ is principal, the line bundle $L_E$ 
is trivial on $\tilde U$, and we may choose a metric $h$ on $L_E$ 
such that $c_1(L_E,h) = 0$.  

We assume from now on that the blow-up is non-trivial, i.e. that 
$\Cal I$ has support of codimension at least 1 in $U$ and is not 
principal on $U$.  In this case  $r>1$ and the   
generators $f_1, .... , f_r$ of 
$\Cal I$ give an embedding  
     $$\iota_f : \tilde U \hookrightarrow U \times \Bbb P^{r-1},$$ 
as described in section III.  Let $[\xi_1: ... : \xi_r]$ be 
homogeneous coordinates for $\Bbb P^{r-1}$.  The blow-up $\tilde U$ 
is covered by open sets 
     $$\tilde U_i  = \tilde U \cap \{\xi_i \neq 0 \}$$ 
on which $L_E$ is trivial.  Transition functions for $L_E$ on the 
intersections $\tilde U_i \cap \tilde U_j$ are the functions 
$g_{ij} = \frac {\xi_i} {\xi_j}$.  
To distinguish between a generating function $f_i$ on $U$ 
and its pullback to $\tilde U \subset U \times \Bbb P^{r-1}$, 
we will let 
     $$\tilde f_i = \pi^* f_i.$$
The exceptional divisor $E$ is given on $\tilde U_i$ by $\tilde f_i = 0$. 
The collection of functions $\tilde f_i$ on the sets $\tilde U_i$ determines 
a section of $L_E$ over $\tilde U$, vanishing exactly on $E$.

\proclaim{Lemma VII.1}  Let $U$ be an open set in $\Bbb C^n$ and let 
$\pi : \tilde U \rightarrow U$ be the blow-up of $U$ along a coherent 
sheaf of ideals $\Cal I$ which is generated by global sections $f_1, ... , f_r$ 
on $U$.  Suppose that the blow-up is non-trivial and 
that $\tilde U$ is smooth.  
Let  $E$ be the exceptional divisor and 
$L_E$ the associated line bundle on $\tilde U$.  Then the embedding 
$\iota_f : \tilde U \hookrightarrow U \times \Bbb P^{r-1}$ 
induces a metric $h$ on $L_E$ whose Chern 
form $c_1 (L_E,h)$  is negative semi-definite on $\tilde U$,  
negative definite on the fibres of the map $E \rightarrow C = V(\Cal I)$,  
and given on  $\tilde U -  E$ by 
     $$c_1 (L_E,h) = 
       \pi^* ( - \rn2pi \ppbar \log \sum_{j=1}^r \absq { f_j (z)} ).$$
\endproclaim

\demo{Proof}  We will construct the Chern form by pullback, 
without using an explicit formula for the metric $h$.  For an 
explicit local formula for $h$, see the remark following this 
proof.  

Recall that the exceptional line bundle $L_E$ on $\tilde U$ 
is the pullback of the universal bundle $\Cal O_{\Bbb P^{r-1}}(-1)$.  
The Fubini-Study form $\omega_{\text{Fub-St}}$ on $\Bbb P^{r-1}$ 
gives a Chern form for $\Cal O_{\Bbb P^{r-1}}(1)$ and $- \omega_{\text{Fub-St}}$ 
gives a Chern form for $\Cal O_{\Bbb P^{r-1}}(-1)$. 
Pulling back to $\tilde U$, we obtain a Chern form for $L_E$ (with 
respect to an induced metric $h$) given 
by 
     $$c_1 (L_E,h)= \iota_f^* \sigma_2^* ( -\omega_{\text{Fub-St}}), $$
where $\sigma_2$ is the second projection map $\sigma_2 :  U \times \Bbb P^{r-1}
\rightarrow \Bbb P^{r-1}$ and $\iota_f$ is the inclusion map 
$\iota_f : \tilde U \hookrightarrow U \times \Bbb P^{r-1}$.  
The negativity properties of $c_1 (L_E,h)$ stated in the lemma 
follow directly from the 
fact that $ \omega_{\text{Fub-St}}$ is positive on $\Bbb P^{r-1}$.

Now recall the formula for the Fubini-Study form on projective space.  If $\xi_1, ...,
\xi_r$ are homogeneous coordinates on $\Bbb P^{r-1}$,  then $w_{ij} = { {\xi_j} \over
{\xi_i}}$ for $j \neq i$ are nonhomogeneous coordinates on $U_i = \{\xi_i \neq 0 \}$. 
The Fubini-Study form $\omega_{\text{Fub-St}}$ is given on $U_i$ by
     $$\align \omega_{\text{Fub-St}} 
     &= \rn2pi \ppbar \log 
     ( 1 + \sum_{j \neq i} \absq {w_{ij}} ) \\
     &= \rn2pi \ppbar \log 
     ( 1 + \sum_{j \neq i} 
       \left \vert \frac {\xi_j} {\xi_i} \right \vert^2 ) . 
     \endalign$$
We continue to use the notation $\tilde f_i = \pi^* f_i$ to 
distinguish between the function $f_i$ on $U$ and its pullback to $\tilde U$.  
On $ \tilde U_i = \tilde U \cap U_i  $ we have  $\frac {\xi_j} {\xi_i}
= \frac {\tilde f_j} {\tilde f_i}$ which gives  us 
     $$\align c_1 (L_E,h) 
     &=  - \rn2pi  \ppbar \log 
     ( 1 + \sum_{j \neq i} \left \vert { {\tilde f_j} \over {\tilde f_i}}
     \right \vert^2 )\\
     &=  - \rn2pi  \ppbar \log 
      {{\sum_{j=1}^r \absq {\tilde f_j}} \over {\absq {\tilde f_i}}} . \endalign$$
On $\tilde U_i - \tilde U_i \cap E$ we have $\tilde f_i(z) \neq 0$ so 
     $$\align c_1 (L_E,h) 
     &=    - \rn2pi  \ppbar 
     ( \log  \sum_{j=1}^r \absq {\tilde f_j}  
     - \log \absq {\tilde f_i} )  \\
     &=    - \rn2pi  \ppbar \log \sum_{j=1}^r \absq { \tilde f_j (z)} \\
     &= \pi^* ( - \rn2pi \ppbar \log 
     \sum_{j=1}^r \absq {  f_j (z)} ) .
     \endalign$$
This formula is independent of $i$, so is valid on all of $\tilde U - E$. 
\qed
\enddemo 

\remark{Remark}  Local defining functions  for the metric $h$ on $L_E$ 
induced from the embedding $\tilde U \hookrightarrow U \times \Bbb P^{r-1}$ 
may also be given explicitly.  
Let $s$ be the section of $L_E$ given on $\tilde U_i$ by $\tilde f_i = 0$.  
The norm of $s$ under the metric $h$ is given by 
     $$\normsq s = \sum_{j=1}^r \absq {\tilde f_j}.$$
The metric $h$ is described locally by positive $C^\infty$ functions $h_i$ 
on $\tilde U_i$ satisfying  
     $$\normsq s = \absq {\tilde f_i} h_i.$$  
Thus 
     $$h_i = \frac {\sum_{j=1}^r \absq {\tilde f_j}} 
       {\absq { \tilde f_i}} .$$
\endremark

\head Global Chern Forms for Blow-ups \endhead

\proclaim{Proposition VII.2}  Let $\pi : \tilde M \rightarrow M$ 
be the blow-up of a compact complex manifold $M$ along a coherent 
sheaf of ideals $\Cal I$ such that $\tilde M$ is smooth.  
Let $E$ be the exceptional divisor 
and $L_E$ the associated line bundle.  

Then there is a metric $h$ on $L_E$ whose 
Chern form $c_1 (L_E,h)$ on $\tilde M$  is negative 
definite along the fibres of the map $E \rightarrow C$ and is 
given on $\tilde M - E$ by 
     $$c_1(L_E,h) = \pi^* ( - \rn2pi \ppbar \log F ),$$
where $F$ is a global $C^\infty$ function on $M$, vanishing 
on the support of $\Cal I$.  Furthermore, $F$ may be constructed 
to be of the form 
     $$F = \prod_\alpha F_\alpha^{\rho_\alpha},$$
where $\{\rho_\alpha\}$ is a $C^\infty$ partition of unity subordinate to an 
open cover $\{ U_\alpha \}$ of $M$,     $F_\alpha$ is a function on 
$U_\alpha$ of the form 
     $$F_\alpha = \sum_{j=1}^r \absq {f_j},$$
and $f_1, ... , f_r$ are local holomorphic generators of the 
coherent ideal sheaf $\Cal I$ on $U_\alpha$. 
\endproclaim

\demo{Proof}  Let $\{U_\alpha \}$  be a finite
open cover of $M$ 
by open sets small enough that $\Cal I$ is generated 
by global sections on each $U_\alpha$.    
If the support of $\Cal I$ 
does not intersect some $U_\alpha$ or if $\Cal I$ is generated by a single generator 
on $U_\alpha$, then $L_E$ is trivial on 
the set $\tilde U_\alpha = \pi^{-1} (U_\alpha)$.  In this case we may 
choose $F_\alpha$ to be a constant and the local Chern form will be 0.  Otherwise, 
in the nontrivial case, suppose that 
$f_1, ... , f_r$ are local generating functions for $\Cal I$ 
on $U_\alpha$ and let 
       $$F_\alpha = \sum_{j=1}^r \absq {f_j} 
       \qquad \text{and} \qquad \tilde F_\alpha = \pi^* F_\alpha.$$   
By Lemma VII.1, there is a local $C^\infty$ metric $h_\alpha$ for $L_E$ 
on $\tilde U_\alpha $ which is negative definite 
on the fibres of the map $E \rightarrow C$ and is 
given on  $\tilde U_\alpha - \tilde U_\alpha \cap E$ 
by 
     $$c_1 (L_E, h_\alpha) 
     = \pi^* ( - \rn2pi \ppbar \log  F_\alpha )
       = - \rn2pi \ppbar \log \tilde F_\alpha . $$

Now choose a $C^\infty$ partition of unity $\{ \rho_\alpha\}$ 
subordinate to $\{ U_\alpha \}$ and let $\trho$ be the pullback of $\rho_\alpha$ to
$\tilde M$.  Then 
$\{ \trho \}$ is a partition of unity on $\tilde M$ subordinate to the open
sets $\{ \tilde U_\alpha \}$.  Note that each function $\trho$ is constant along the fibres
of the map $E \rightarrow C$.  

We define a global $C^\infty$ metric for $L_E$ as follows.  
For any section $s$  of $L_E$, let $\mid \mid s \mid \mid^2_\alpha$ be 
the norm-squared of 
$s$ with respect to the metric $h_\alpha$ on $\tilde U_\alpha$ and let 
     $$\normsq s = \prod_\alpha \mid \mid s \mid \mid_\alpha^{2 {\tilde \rho}_\alpha}.$$
If $\{V_i \}$ is a cover of $\tilde M$ by open sets on which $L_E$ is 
trivial, and $h_{\alpha i}$ is the positive $C^\infty$ function representing $h_\alpha$ on 
$V_i$, then the positive $C^\infty$ function for $h$ on $V_i$ is 
     $$h_i = \prod_\alpha h_{\alpha i}^{{\tilde \rho}_\alpha}.$$
If $s$ is given on $V_i$ by the holomorphic function $s_i$, then 
on $U_\alpha \cap V_i$ we have 
     $$\mid \mid s \mid \mid^2_\alpha = \absq {s_i} h_{\alpha i}$$ 
and on $V_i$, 
     $$\normsq s = \absq {s_i} h_i.$$

The global form $c_1(L_E,h)$ associated with this metric 
preserves the property of the local forms 
of being negative definite on the fibres of the map $E \rightarrow C$
because the partition of unity functions $\{ \trho \}$ are  constant on fibres of the map $E
\rightarrow C$.  On $V_i$, this Chern form is given by 
     $$\align c_1 (L_E, h) 
     &= - \rn2pi \ppbar \log h_i \\
       &= - \rn2pi \ppbar \log \prod_\alpha h_{\alpha i}^{{\tilde \rho}_\alpha} \\
     &= - \rn2pi \sum_\alpha \ppbar  {\tilde \rho}_\alpha \log h_{\alpha i}. \endalign$$

Let $s$ be a global holomorphic section of $L_E$ on $\tilde M$ whose 
associated divisor is $E$.  Such a section always exists - just choose 
local holomorphic defining equations of $E$ to determine $s$ locally.  
For example, on $\tilde U_{\alpha i} = \tilde U_\alpha \cap \{ \xi_i \neq 0 \} 
\subset U_\alpha \times \Bbb P^{r-1}$, 
take  $ s_{\alpha i} = \tilde f_i = \pi^* f_i$, where $f_1, ... , f_r$ are 
local holomorphic generators of $\Cal I$ on $U_\alpha$.  
Then 
     $$\mid \mid s \mid \mid^2_\alpha = \sum_{j=1}^r \absq {\tilde f_j} 
     = \tilde F_\alpha  $$
and 
     $$\normsq s = \prod_\alpha \tilde F_\alpha^{{\tilde \rho}_\alpha}
       = \pi^* ( 
       \prod_\alpha  F_\alpha^{{ \rho}_\alpha} ).$$
Thus the Chern form $c_1(L_E,h)$ is given on $\tilde M - E$ by 
     $$c_1(L_E,h) = - \rn2pi \ppbar \log \normsq s 
       = \pi^* ( - \rn2pi \ppbar \log F ),$$ 
where 
      $$F = \prod_\alpha F_\alpha^{{\rho}_\alpha}. \qed$$
\enddemo

\head VIII.  Construction of Saper-Type Metrics \endhead

Let $X$ be a singular subvariety of a compact K\"ahler manifold 
$M$ and let $\xsing$ be the singular locus of $X$.  We will construct 
Saper-type metrics on $M - \xsing$, first locally, then globally  
using a $C^\infty$ partition of unity on $M$.  These metrics are 
complete K\"ahler metrics on $M - \xsing$ which grow less rapidly 
than Poincar\'e metrics near the singular locus.  More details 
on the growth rate of Saper-type metrics and their relationship to 
intersection cohomology may be found in [GM], [Sa1], and [Sa2].  

      We also construct a non-complete K\"ahler metric 
on $M - \xsing$ with the property that the completion of $X - \xsing$ 
with respect to this metric is a desingularization of $X$.  We call 
this metric a \lq \lq desingularizing metric" for $X$.  

     The constructions of both metrics are based on 
resolution of singularities using a single coherent ideal 
sheaf $\Cal I$ on $M$ (see Corollary VI.3) and the 
explicit formula for a Chern form for the blow-up of $M$ 
along $\Cal I$ given in Proposition VII.2

\head Local Construction of Metrics \endhead

Before constructing Saper-type metrics, we will describe a 
K\"ahler metric for a local blow-up.   
                                   
Let $U$ be an open set in $\Bbb C^n$ and let 
$\pi : \tilde U \rightarrow U$ be the blow-up of  $U$ 
along a coherent sheaf of ideals $\Cal I$ such that 
$\tilde U$ is smooth.  
Let $E$ be the exceptional divisor of $\pi$.
Assume that $U$ 
is small enough that $\Cal I$ is generated by global sections on $U$ 
and let  
     $$\iota_f : \tilde U \hookrightarrow \UP r$$
be the embedding associated with a collection of generators $f$.  Let 
$\sigma_1$ and $\sigma_2$ be the projection maps 
$$\CD   U \times \Bbb P^{r-1} @> \sigma_2 >>
     \Bbb P^{r-1} \\
     @V \sigma_1 VV @. \\
     U  @. {}
     \endCD$$
 
Suppose that $\omega$ is the K\"ahler
form of a K\"ahler metric on $U$ and let $\FubS$ be the K\"ahler form of the Fubini-Study
metric on $\Bbb P^{r-1}$.  

\proclaim{Lemma VIII.1} 
The embedding $\tilde U \hookrightarrow U \times \Bbb P^{r-1}$ 
induces a K\"ahler metric on $\tilde U$ whose K\"ahler form is 
     $$\omega^\prime = \pi^* \omega - c_1 ( L_E,h),$$
where $c_1 (L_E,h)$ is a Chern form of the line bundle $L_E$ (with 
respect to a metric $h$) of the type described in Lemma VII.1. 
\endproclaim

\demo{Proof}  The K\"ahler form on $\tilde U$ given by the restriction of 
the product metric on $U \times \Bbb P^{r-1}$ is 
     $$\align  \omega^\prime  
     &= \iota_f^* ( \sigma_1^* \omega + \sigma_2^* \FubS )\\
     &  = \pi^* \omega + \iota_f^* \sigma_2^* \FubS \\
     &= \pi^* \omega - c_1 (L_E,h). \qed	\endalign$$
\enddemo

Recall that the Chern form of Lemma VII.1 was given on the set 
$\tilde U - E$ by 
     $$c_1 (L_E,h) = \pi^* ( - \rn2pi \ppbar 
       \log \sum_{j=1}^r \absq {f_j} ),$$
where $f_1, ... , f_r$ were holomorphic generators for $\Cal I$ 
on $U$.  Thus the (1,1)-form 
     $$\tilde \omega = \omega + \rn2pi 
     \ppbar \log \sum_{j=1}^r \absq {f_j}$$ 
on $U - V( \Cal I)$ determines a K\"ahler metric on $U - V(\Cal I).$ 
This K\"ahler metric is essentially the local model of 
our desingularizing metric. 

The function 
$F = \sum_{j=1}^r \absq {f_j}$ 
can also be used to construct a Saper-type metric on $U - V(\Cal I)$.  
We are particularly interested in the case of a 
coherent sheaf of ideals  $\Cal I$ which 
determines a resolution of singularities of a singular variety 
and which is supported on the singular locus of the variety.   
Theorems VIII.2 and VIII.3 describe  local and global 
constructions, respectively, of Saper and desingularizing  
metrics for a singular variety.  
The main differences between the two theorems are that 
we must patch with a $C^\infty$ partition of unity in the 
global case, and that our desingularizing metric may require 
a multiple of the original metric in that case.    

\proclaim{Theorem VIII.2.  Local Metrics}
Let $X$ be a singular subvariety of a compact K\"ahler manifold $M$ with singular locus
$\xsing$.  Let $\omega$ be the K\"ahler (1,1)-form of a K\"ahler metric on $M$.  Let $p$
be any point in $\xsing$.  Then there exists a neighborhood $U$ of $p$ and a $C^\infty$
function $F$ on $U$, vanishing on $U \cap \xsing$,  such that 
\roster
     \item"{i.}"  the (1,1)-form 
     $$\tilde \omega = \omega + \rn2pi \ppbar \log F$$
is the K\"ahler form of an incomplete metric on $U - U \cap \xsing$ which determines  a local embedded resolution of singularities and 
   \item"{ii.}"  the (1,1)-form 
     $$\omega_S = l \omega - \rn2pi \ppbar \log ( \log F)^2$$ 
on $U - U \cap \xsing$ is the K\"ahler form of a modified Saper metric on $U - U \cap
\xsing$ (in the sense of [GM]) for $l$ a large enough positive integer.
\endroster  
Furthermore, the function $F$ may be constructed to be of the form 
     $$F = \sum_{j=1}^r \absq {f_j},$$
where $f_1, ... , f_r$ are holomorphic functions on $U$ which are 
local generators of a coherent ideal sheaf $\Cal I$ on $M$, such that  
blowing up $M$ along 
$\Cal I$ desingularizes $X$, 
$\Cal I$ is supported on $\xsing$, 
and the exceptional 
divisor of the blow-up along $\Cal I$ has normal crossings with itself and 
with the desingularization of $X$.      
\endproclaim

\demo{Proof}  Part (i) is a consequence of Lemma VIII.1 and Lemma 
VII.1.  Part (ii) and its global version follow from Theorem 9.2.1 of [GM], in 
which we also give estimates of the rate of growth of $\omega_S$.  
The idea behind the proof of that theorem is that we can 
decompose the term involving $F$ in our Saper-type metric as 
     $$ - \rn2pi \ppbar \log ( \log F)^2 
     = \frac { \sqrt {-1} } {2 \pi} 
       \left( \frac 1 {\abs {\log F}} \ppbar \log F 
     + \frac { \partial F \wedge \overline \partial F} 
       { \absq F ( \log F)^2} \right).$$
The first term gives positivity of $\omega_S$, since it 
is a multiple of the negative of the Chern form $c_1 (L_E,h)$.  
The second term is similar to the Poincar\'e metric on the punctured 
disc and becomes unbounded near the singular locus.  
\qed
\enddemo

\head Global Construction of Metrics \endhead

To construct global metrics we patch together our local metrics using $C^\infty$
partitions of unity on $M$.  As described in section VII, this patching does not affect
negativity of our Chern forms along fibres of the maps $E \rightarrow C$ from  the
exceptional divisors to their corresponding centers.  However the Chern forms may not
remain negative semidefinite elsewhere, so that it may be necessary to introduce
multiples of the original metric. 

\proclaim{Theorem VIII.3.  Global Metrics}
Let $X$ be a singular subvariety of a compact K\"ahler manifold $M$ with singular locus
$\xsing$.  Let $\omega$ be the K\"ahler (1,1)-form of a K\"ahler metric on $M$.  There 
exists a global $C^\infty$ function $F$ on $M$, vanishing exactly on $\xsing$, such that 
for $l$ a large enough positive integer 
\roster
   \item"{i.}"  the (1,1)-form 
     $$\tilde \omega = l \omega + \rn2pi \ppbar \log F$$
is the K\"ahler form of an incomplete K\"ahler metric on $M - \xsing$ 
which is a desingularizing metric for $X$ (i.e. the completion 
of $X - \xsing$ with respect to $\tilde \omega$ is nonsingular), and 
   \item"{ii.}"  the (1,1)-form 
          $$\omega_S = l \omega - \rn2pi \ppbar \log ( \log F)^2$$ 
on $M - \xsing$ is the K\"ahler form of a complete K\"ahler 
modified Saper metric (in the sense of [GM]). 
\endroster
Furthermore, the function $F$ may be constructed 
to be of the form 
     $$F = \prod_\alpha F_\alpha^{\rho_\alpha},$$ 
where $\{ \rho_\alpha \}$ is a $C^\infty$ partition of 
unity subordinate to an open cover $\{U_\alpha \}$ of $M$, 
$F_\alpha$ is a function on $U_\alpha$ of the form 
     $$F_\alpha = \sum_{j=1}^r \absq {f_j},$$ 
and $f_1, ... , f_r$ are holomorphic functions on $U_\alpha$, 
vanishing exactly on $\xsing \cap U_\alpha$.  
More specifically, $f_1, ... , f_r$ are local holomorphic 
generators of a coherent ideal sheaf $\Cal I$ on $M$ 
such that blowing up $M$ along 
$\Cal I$ desingularizes $X$, 
$\Cal I$ is supported on $\xsing$, 
and the exceptional 
divisor of the blow-up along $\Cal I$ has normal crossings with itself and 
with the desingularization of $X$.     
\endproclaim

\demo{Proof}  Part (i) follows from Proposition VII.2.  Part (ii) follows from 
our description of $F$ in section VII.2 and Theorem 9.2.1 of [GM].
\qed
\enddemo

\head IX.  Example \endhead

\subhead  The cuspidal cubic  \endsubhead 
Let $M = \Bbb P^2$ and let $X$ be the cuspidal cubic given in homogeneous 
coordinates by $\xi_0 \xi_2^2 - \xi_1^3 = 0$.  In local coordinates $x,y$ 
in a neighborhood $U \cong \Bbb C^2$ of the singular point,  $X$ is 
given  by 
     $$y^2 - x^3 = 0.$$ 
The singularity may be resolved by three blow-ups of points, in such a 
way that  the components of the total exceptional divisor 
have normal crossings with each other and with the desingularization of $X$. 
We will show that these three blow-ups are equivalent to a single 
blow-up along the ideal sheaf given locally by 
     $$\Cal I = (x,y)(x^2,y)(x^3,x^2y,y).$$

\subhead  First blow-up $\pi_1$  \endsubhead
The center $C_1$ for the first blow-up is the point $x=y=0$ and its ideal 
is $\Cal I_{C_1} = (x,y)$.  
The blow-up $U_1 = \pi_1^{-1} (U)$ may 
be covered by two coordinate charts, which we will call the $x$- and 
$y$-coordinate charts, according to whether 
the chart is a complement in $U_1$ of the strict transform of $x=0$ or 
$y=0$.    
(The exceptional divisor 
is given by the vanishing of the $x$-coordinate in the $x$-chart 
and the $y$-coordinate in the $y$ chart.)
On the $x$-coordinate 
chart, $\pi_1$ is given by 
     $$\pi_1(x_1, y_1) = (x_1, x_1 y_1) = (x,y)$$ 
and the exceptional divisor $E_1$ is given by $x_1 = 0$. 
The inverse image $\pi_1^{-1} (X)$ is given by 
$x_1^2 y_1^2 - x_1^3 = 0$.  
The strict transform $X_1$ of $X$ is obtained from the inverse image 
by removing all copies 
of $E_1$,  i.e. by dividing by the highest possible power of $x_1$, 
which gives 
     $$y_1^2 - x_1 = 0.$$
Although $X_1$ is smooth, it does not have normal crossings with 
the divisor $E_1$ at the point $x_1 = y_1 = 0$, so we must blow up 
again at this point.  Before doing so, we note that in the $y$-coordinate 
chart, the strict transform $X_1$ is smooth and has normal 
crossings with $E_1$, so there is no need to blow up further at any 
points in that chart.  

\subhead Second blow-up $\pi_2$ \endsubhead
The center $C_2$ for the second blow-up is the point $x_1 = y_1 = 0$ 
in the $x$-coordinate chart of $U_1$, and its ideal is 
$\Cal I_{C_2} = (x_1,y_1)$.  
In the $x$-coordinate chart of $\pi_2$ we have normal crossings,  
so there is no need to blow up further at any points in that  
chart.   
In local coordinates $(x_2, y_2)$ for the $y$-coordinate chart of $\pi_2$, we have 
     $$\pi_2(x_2, y_2) = (x_2 y_2, y_2) = (x_1, y_1)$$ 
and $\Cal I_{E_2} = (y_2)$.  
The strict transform $X_2$ of $X_1$ is given by 
     $$y_2 - x_2 = 0$$ 
and the strict transform $\tilde E_1$ of $E_1$ by $x_2 = 0$.  
The total 
exceptional divisor of the first two blow-ups, which is the union of 
$E_2$ and $ \tilde E_1$, does not 
have normal crossings with $X_2$ so we blow up again.  

\subhead Third blow-up $\pi_3$ \endsubhead
The center $C_3$ for the third blow-up is the point $x_2 = y_2 = 0$ 
with ideal $\Cal I_{C_2} = (x_2,y_2)$.  
After this third blow-up, the strict transform of $X$ and all 
three components of the total exceptional divisor have normal 
crossings.

\subhead  Construction of $\Cal I$  \endsubhead
We will construct $\Cal I$ as a product 
$\Cal I = \Cal I_1 \Cal I_2 \Cal I_3$ 
of ideals corresponding to the centers of the blow-ups.  
We begin by 
choosing $\Cal I_1 = 
\Cal I_{C_1} = (x,y)$.

To obtain $\Cal I_2$, we start with $\Cal I_{C_2} $ and multiply by a 
high enough power of $\Cal I_{E_1} $ such that taking the direct image 
under $\pi_1$ and then the inverse image does not change the ideal.  
We define $\Cal I_2$ to be the direct image of the resulting product 
under the map $\pi_1$. 
 
Locally, in the $x$-coordinate chart of $\pi_1$, 
$\Cal I_{C_2}$ is given by $(x_1, y_1)$ and $\Cal I_{E_1}$ 
by $(x_1)$, 
where $x_1 = x$ and $y_1 = \frac y x$. 
Thus $\Cal I_{C_2}$ is not the inverse image of an ideal 
sheaf, but $\Cal I_{C_2} \Cal I_{E_1}$ is, since 
     $$\pi_1^{-1} (x^2,y) = \Cal I_{C_2} \Cal I_{E_1}.$$  
 The direct image 
$\pi_{1*} ( \Cal I_{C_2} \Cal I_{E_1})$ is the largest ideal sheaf whose 
inverse image is contained in $\Cal I_{C_2} \Cal I_{E_1}$, so 
$\pi_{1*} ( \Cal I_{C_2} \Cal I_{E_1}) $ contains  $(x^2,y)$.  It is 
easily checked that $x^2$ and $y$ generate 
$\pi_{1*} ( \Cal I_{C_2} \Cal I_{E_1}) $, since
they 
are the only monomials whose 
pullbacks are sections of $ \Cal I_{C_2} \Cal I_{E_1}$.  Thus 
     $$\Cal I_2 = \pi_{1*}(\Cal I_{C_2} \Cal I_{E_1}) = (x^2,y).$$ 
      
Similarly, to obtain $\Cal I_3$ we start with $\Cal I_{C_3}$, 
given locally by  $ (x_2,y_2)$, and 
recall that $x_2 = \frac {x_1}{y_1}$ and $y_2 = y_1 $.  
Hence $\Cal I_{C_3} \Cal I_{E_2}$ is the inverse 
image of an ideal sheaf $\Cal J$ given locally on $U_1$ by $(x_1, y_1^2)$, 
and $\Cal J \Cal I_{E_1}^2$ is the inverse image of the ideal sheaf 
$(x^3,y^2)$.  Since $\pi_2^{-1} (\Cal I_{E_1}) = \Cal I_{\tilde E_1 } 
\Cal I_{E_2}$, it follows that 
     $$\pi_2^{-1} \pi_1^{-1} (x^3,y^2)
= \Cal I_{C_3} \Cal I_{\tilde E_1}^2 \Cal I_{E_2}^3.$$ 
In local coordinates, 
$\pi_2^{-1} \pi_1^{-1} (x^3,y^2) =(x_2,y_2)(x_2^2)(y_2^3)$.  
We define 
$\Cal I_3$ to be the direct image 
$\pi_{1*} \pi_{2*} (\Cal I_{C_3} \Cal I_{\tilde {E_1}}^2 \Cal I_{E_2}^3)$,  
and note that $\Cal I_3 $ contains $ (x^3,y^2)$,
since $\Cal I_3$ is the largest ideal sheaf whose inverse 
image is contained in $\Cal I_{C_3} \Cal I_{\tilde {E_1}}^2 \Cal I_{E_2}^3$.    
To find any remaining generators of $\Cal I_3$, we test monomials 
not generated by $x^3$ or $y^2$ to see which pull back to sections of 
$\Cal I_{C_3} \Cal I_{\tilde {E_1}}^2 \Cal I_{E_2}^3$.  
It is easily checked that $x$, $y$, $x^2$, and $xy$ are not in $\Cal I_3$, 
 but $x^2y$ is in $\Cal I_3$ since $x^2 y = x_1^3y_1 = x_2^3y_2^4$.  Thus 
     $$\Cal I_3 =  \pi_{1*} \pi_{2*} (
     \Cal I_{C_3} \Cal I_{\tilde {E_1}}^2 \Cal I_{E_2}^3) 
     = (x^3, x^2y , y^2).$$

We define the ideal $\Cal I$ to be the product of $\Cal I_1$, $\Cal I_2$, and 
$\Cal I_3$ 
      $$\Cal I = (x,y) (x^2, y) (x^3, x^2y, y^2).$$
Blowing up along $\Cal I$ is equivalent to blowing up sequentially along 
the centers $C_1$, $C_2$, and $C_3$.  

    The method used in this example may be generalized to any locally toric  
complex analytic variety.  Details will be given elsewhere.


\Refs

\widestnumber\key{ABS1}

\ref \key BM \by E. Bierstone and P. Milman \pages 207--302   
   \paper Canonical desingularization in characteristic zero by 
   blowing up the maximum strata of a local invariant
   \yr1997 \vol 128 \jour Invent. Math. \endref

\ref \key C \by J. Cheeger 
   \paper On the Hodge Theory of Riemannian Pseudomanifolds
   \yr 1980  \vol 36  \jour Proc. Symp. Pure Math., American 
   Math. Soc.  \pages 91--146
   \endref

\ref \key CGM  \by J. Cheeger, M. Goresky, and R. MacPherson
   \paper $L^2$-Cohomology and Intersection Homology of 
   Singular Algebraic Varieties  
   \inbook Seminar on Differential Geometry 
   \bookinfo Annals of Mathematics Studies 102
   \publ Princeton University Press \publaddr Princeton, NJ 
   \yr 1982  \pages 303--340 
   \endref

\ref \key F \by G. Fischer 
   \book Complex Analytic Geometry 
   \bookinfo Lecture Notes in Math.  538
   \publ Springer-Verlag 
   \publaddr Berlin Heidelberg 
   \yr 1976   \endref

\ref \key Ful \by W. Fulton
   \book Intersection Theory
   \bookinfo Ergebnisse der Mathematik und ihrer 
     Grenzgebiete; 3. Folge, Bd 2
   \publ Springer-Verlag 
   \publaddr Berlin Heidelberg 
   \yr 1984   \endref 

\ref \key GH \by P. Griffiths and J. Harris
   \book Principles of Algebraic Geometry
   \publ Wiley-Interscience \publaddr New York \yr1978 \endref

\ref \key GM \by C. Grant and P. Milman  \pages 61--156
   \paper Metrics for Singular Analytic Spaces
   \yr 1995  \vol 168  \jour  Pac. J. Math.    \endref

\ref \key GoM \by M. Goresky and R. MacPherson \pages 77--129
   \paper Intersection Homology II
   \yr1983 \vol 71 \jour Invent. Math. \endref

\ref \key GrR1 \by H. Grauert and R. Remmert  
   \book Coherent Analytic Sheaves  
   \bookinfo Grundlehren der mathematischen Wissenschaften 265
   \publ Springer-Verlag \publaddr Berlin Heidelberg 
   \yr 1984  \endref

\ref \key GrR2 \by H. Grauert and R. Remmert  
   \book Theory of Stein Spaces  
  \bookinfo Grundlehren der mathematischen Wissenschaften 236
   \publ Springer-Verlag \publaddr New York 
   \yr 1979  \endref


\ref \key GuR \by R. Gunning and H. Rossi 
   \book Analytic Functions of Several Complex Variables
   \publ Prentice-Hall Inc. \publaddr Englewood Cliffs, NJ
   \yr 1965  \endref

\ref \key Ha1 \by R. Hartshorne 
   \book Algebraic Geometry  \bookinfo Graduate Texts in 
   Mathematics 52 \publ Springer-Verlag
   \publaddr New York \yr 1977
   \endref

\ref \key Ha2 \by R. Hartshorne 
   \book  Ample Subvarieties of Algebraic Varieties
   \bookinfo Lecture Notes in Math. 156
   \publ Springer-Verlag \publaddr Heidelberg \yr 1970
   \endref

\ref \key Hi \by H. Hironaka 
   \paper Resolution of singularities of an algebraic
   variety over a field of characteristic zero: I, II
   \yr1964 \vol 79 \jour Ann. Math.  \pages 109--326 \endref

\ref \key H\"o \by L. H\"ormander
     \book An Introduction to Complex Analysis in Several Variables 
     \publ North-Holland \publaddr New York \yr 1973 
     \endref

\ref \key HR \by H. Hironaka and H. Rossi  
   \paper On the Equivalence of Imbeddings of 
   Exceptional Complex Spaces
   \yr 1964  \vol 156  \jour Math. Annalen \pages 313--333
   \endref

\ref \key I \by S. Iitaka
   \book  Algebraic Geometry \bookinfo Graduate Texts in 
   Mathematics 76  \publ Springer-Verlag \publaddr 
   New York \yr 1982
   \endref

\ref \key K \by E. Kunz \book Introduction to Commutative 
   Algebra and Algebraic Geometry \publ Birkh\"auser 
   \publaddr Boston \yr 1985
   \endref

\ref \key M \by D. Mumford 
   \book Algebraic Geometry I Complex Projective Varieties
   \bookinfo Grundlehren der mathematischen Wissenschaften 221 
   \publ Springer-Verlag \publaddr Berlin Heidelberg 
   \yr 1976
   \endref

\ref \key Ma \by H. Matsumura   
   \book Commutative Algebra 
   \publ W. A. Benjamin Co. \publaddr New York  
   \yr 1970
   \endref

\ref \key Sa1 \by L. Saper \pages 207--255
   \paper $L_2$-cohomology and intersection homology of certain algebraic 
   varieties with isolated singularities 
   \yr1985 \vol 82 \jour Invent. Math. \endref

\ref \key Sa2 \by L. Saper \paper $L_2$-cohomology of K\"ahler
   varieties with isolated singularities 
  \jour J. Diff. Geom. \vol 36  \yr1992 \pages 89--161
   \endref

\ref \key Sh \by I. Shafarevich 
   \book Basic Algebraic Geometry  2 
   \publ Springer-Verlag
   \publaddr Berlin Heidelberg
   \yr 1994 \endref

\ref \key Sp \by M. Spivakovsky \pages 107--156
   \paper Valuations in Function Fields of Surfaces
   \yr 1990  \vol 112  \jour Am. J. Math. \endref

\ref \key W \by R. O. Wells 
   \book Differential Analysis on Complex Manifolds
   \bookinfo Graduate Texts in Mathematics 65
   \publ Springer-Verlag \publaddr New York  \yr1980 \endref



\endRefs
\enddocument